\documentclass[11pt,a4paper]{article}
\usepackage{amsfonts,amsmath,amssymb,accents,bm,amsthm}
\usepackage{verbatim}

\usepackage{hyperref}
\usepackage[normalem]{ulem}

\usepackage{tikz} 
\usetikzlibrary{math,calc,arrows,positioning,fit,petri,external}
\newcommand{\inputtikz}[1]{\includegraphics{figures/#1}}

\usepackage{graphicx}

\newcommand{\dis}{\displaystyle}

\newcommand{\noi}{\noindent}
\newcommand{\halmos}{\rule{1ex}{1.4ex}}
\newcommand{\QED}{\nopagebreak{\hspace*{\fill}$\halmos$\medskip}}
\newcommand{\med}{\medskip}
\newcommand{\quand}{\quad\mbox{and}\quad}

\newtheoremstyle{mythm}
  {}
  {}
  {\itshape}
  {}
  {\bfseries}
  {}
  {.5em}
  {#1 #2 \thmnote{(#3)}}

\theoremstyle{mythm}
\newtheorem{theorem}{Theorem}
\newtheorem{proposition}[theorem]{Proposition}
\newtheorem{lemma}[theorem]{Lemma}

\newcommand{\bt}{\begin{theorem}}
\newcommand{\et}{\end{theorem}}
\newcommand{\bl}{\begin{lemma}}
\newcommand{\el}{\end{lemma}}
\newcommand{\bp}{\begin{proposition}}
\newcommand{\ep}{\end{proposition}}

\newenvironment{Proof}[1][]{\noi\textbf{Proof #1}}{\QED}
\newcommand{\bpro}{\begin{Proof}}
\newcommand{\epro}{\end{Proof}}

\newcommand{\be}{\begin{equation}}
\newcommand{\ee}{\end{equation}}
\newcommand{\ba}{\begin{array}}
\newcommand{\ea}{\end{array}}
\newcommand{\bac}{\begin{array}{r@{\,}c@{\,}l}}
\newcommand{\bc}{\be\begin{array}{r@{\,}c@{\,}l}}
\newcommand{\ec}{\end{array}\ee}

\newcommand{\de}{\delta}

\newcommand{\eps}{\varepsilon}

\newcommand{\sig}{\sigma}

\newcommand{\si}{\ensuremath{\sigma}}

\newcommand{\Ai}{{\cal A}}

\newcommand{\Di}{{\cal D}}

\newcommand{\Fi}{{\cal F}}

\newcommand{\Mi}{{\cal M}}

\newcommand{\Oi}{{\cal O}}

\newcommand{\Si}{{\cal S}}

\newcommand{\Zi}{{\cal Z}}

\newcommand{\N}{{\mathbb N}}

\renewcommand{\P}{{\mathbb P}}

\newcommand{\T}{{\mathbb T}}

\newcommand{\Z}{{\mathbb Z}}

\newcommand{\li}{\langle}
\newcommand{\re}{\rangle}

\newcommand{\desd}{\ensuremath{\Leftrightarrow}}
\newcommand{\volgt}{\ensuremath{\Rightarrow}}
\newcommand{\up}{\uparrow}
\newcommand{\down}{\downarrow}
\newcommand{\sub}{\subset}

\newcommand{\asto}[1]{\underset{{#1}\to\infty}{\longrightarrow}}
\newcommand{\Asto}[1]{\underset{{#1}\to\infty}{\Longrightarrow}}

\newcommand{\ti}{\tilde}

\newcommand{\ov}{\overline}
\newcommand{\un}{\underline}

\newcommand{\rvec}[1]{\accentset{\rightarrow}{#1}}
\newcommand{\lvec}[1]{\accentset{\leftarrow}{#1}}

\newcommand{\pa}{\partial}
\newcommand{\ffrac}[2]{{\textstyle\frac{{#1}}{{#2}}}}
\newcommand{\dif}[1]{\ffrac{\partial}{\partial{#1}}}

\newcommand{\di}{\mathrm{d}}

\newcommand{\ha}{\ffrac{1}{2}}

\newcommand{\var}{{\rm Var}}

\newcommand{\ibf}{\mathbf{i}}
\newcommand{\jbf}{\mathbf{j}}
\newcommand{\vbf}{\mathbf{v}}
\newcommand{\wurz}{\varnothing}

\newcommand{\Ali}{1}
\newcommand{\Bob}{2}
\newcommand{\hX}{\widehat X}
\newcommand{\hDi}{\widehat{\Di\hspace{-0.5pt}}\hspace{0.5pt}}
\newcommand{\hFi}{\widehat\Fi}

\begin{document}

\makeatletter\@addtoreset{equation}{section}
\makeatother\def\theequation{\thesection.\arabic{equation}}

\renewcommand{\labelenumi}{{\rm (\roman{enumi})}}
\renewcommand{\theenumi}{\roman{enumi}}

\title{A min-max random game on a graph that is not a tree}
\author{Natalia~Cardona-Tob\'on\footnote{Departamento de Estad\'istica, Universidad Nacional de Colombia, Carrera 45 No 26-85,  111321 Bogot\'a, Colombia. ncardonat@unal.edu.co}, Anja~Sturm\footnote{Institute for Mathematical Stochastics, Georg-August-Universit\"at G\"ottingen, Goldschmidtstr.~7, 37077 G\"ottingen, Germany. asturm@math.uni-goettingen.de}, Jan~M.~Swart\footnote{The Czech Academy of Sciences, Institute of Information Theory and Automation, Pod vod\'arenskou v\v{e}\v{z}\'i~4, 18200 Praha 8, Czech republic. swart@utia.cas.cz}} 

\date{\today}

\maketitle

\begin{abstract}\noi
We study a random game in which two players in turn play a fixed number of moves. For each move, there are two possible choices. To each possible outcome of the game we assign a winner in an i.i.d.\ fashion with a fixed parameter $p$. In the case where all different game histories lead to different outcomes, a classical result due to Pearl (1980) says that in the limit when the number of moves is large, there is a sharp threshold in the parameter $p$ that separates the regimes in which either player has with high probability a winning strategy. We are interested in a modification of this game where the outcome is determined by the exact sequence of moves played by the first player and by the number of times the second player has played each of the two possible moves. We show that also in this case, there is a sharp threshold in the parameter $p$ that separates the regimes in which either player has with high probability a winning strategy. Since in the modified game, different game histories can lead to the same outcome, the graph associated with the game is no longer a tree which means independence is lost. As a result, the analysis becomes more complicated and open problems remain.
\end{abstract}

\noi
%
%
%

\noi

{\setlength{\parskip}{-2pt}\tableofcontents}

\newpage

\section{Introduction and main results}

\subsection{Main results}\label{S:intro}

Consider a game played by two players, Alice and Bob, who take turns to play $n$ moves each. Alice starts. For each move, each player has two options. The outcome of the game is determined by the exact sequences of moves played by Alice and Bob. We assign a random winner to each of the $2^{2n}$ possible outcomes of the game in an i.i.d.\ fashion. For each given outcome, the probability that Bob is the winner is $p$. For later reference, we give this game the name ${\rm AB}_n(p)$. Since the game is finite, for each possible assignment of winners to outcomes, precisely one of the players has a winning strategy.\footnote{This is a basic result from game theory, that can easily be proved using the inductive formula (\ref{minmax}) below.} Let $P^{\rm AB}_n(p)$ denote the probability that Bob has a winning strategy. Pearl \cite{Pea80} showed that, setting $p^{\rm AB}_{\rm c}:=\ha(3-\sqrt{5})\approx 0.382$, which has the effect that $p^{\rm AB}_{\rm c}:1-p^{\rm AB}_{\rm c}$ is the golden ratio,\footnote{In fact, setting $\beta:=\ha(1+\sqrt{5})$, one has that $p^{\rm AB}_{\rm c}=\beta^{-2}$ and $1-p^{\rm AB}_{\rm c}=\beta^{-1}$.} one has that
\be\label{Pearl}
P^{\rm AB}_n(p)\asto{n}\left\{\ba{ll}
0\quad&\mbox{if }p<p^{\rm AB}_{\rm c},\\
p^{\rm AB}_{\rm c}\quad&\mbox{if }p=p^{\rm AB}_{\rm c},\\
1\quad&\mbox{if }p>p^{\rm AB}_{\rm c}.
\ea\right.
\ee
Note that $p^{\rm AB}_{\rm c}<1/2$, which is due to the fact that Bob has the last move, which gives him an advantage.

Imagine now that we change the game in such a way that for the outcome of the game, the whole history of the moves played by Alice is relevant as before, but all that matters of Bob's moves is the total number of times he plays each of the two possible moves. In this case, there are $2^n(n+1)$ possible outcomes of the game. As before, we assign a winner independently to each possible outcome, where $p$ is the probability that Bob is the winner for a given outcome. We call this game ${\rm Ab}_n(p)$ and let $P^{\rm Ab}_n(p)$ denote the probability that Bob has a winning strategy. Similarly, we let ${\rm aB}_n(p)$ denote the game in which the outcome is determined by how often Alice has played each of the two possible moves and by the whole history of the moves played by Bob. Finally, we let ${\rm ab}_n(p)$ denote the game in which both for Alice and Bob, all that matters for the outcome is how often each of them has played each of the possible moves. We denote the probability that Bob has a winning strategy in these games by $P^{\rm aB}_n(p)$ and $P^{\rm ab}_n(p)$, respectively.

It seems that if for the outcome of a game, all that matters is how often a player has played each of the two possible moves, then this player has less influence on the outcome compared to the game in which the whole history of moves played by this player matters. It is natural to conjecture that this gives such a player a disadvantage. Our first result partially confirms this.

\bp[Less influence on the outcome]
One\label{P:treecomp} has $P^{Ab}_n(p)\leq P^{AB}_n(p)$ and $P^{aB}_n(p)\geq P^{AB}_n(p)$ for all $p\in[0,1]$ and $n\geq 1$.
\ep

We conjecture that similarly $P^{ab}_n(p)\leq P^{aB}_n(p)$ and $P^{ab}_n(p)\geq P^{Ab}_n(p)$, but we have not been able to prove this.

From now on, we mostly focus on the games ${\rm Ab}_n(p)$ and ${\rm aB}_n(p)$. Our next result says that even though these games give Bob and Alice a disadvantage compared to the game ${\rm AB}_n(p)$, for values of $p$ close enough to one or zero, they still have with high probability a winning strategy. Moreover, similarly to (\ref{Pearl}), there is a sharp threshold of the parameter $p$ that separates two regimes where for large $n$ either Alice or Bob have with high probability a winning strategy.

\bt[Sharp threshold]
There\label{T:threshold} exist constants $0<p^{\rm Ab}_{\rm c},p^{\rm aB}_{\rm c}<1$ such that
\be\ba{l}
\dis P^{\rm Ab}_n(p)\asto{n}\left\{\ba{ll}
0\quad&\mbox{if and only if }p<p^{\rm Ab}_{\rm c},\\
1\quad&\mbox{if }p>p^{\rm Ab}_{\rm c},
\ea\right.\\[25pt]
\dis P^{\rm aB}_n(p)\asto{n}\left\{\ba{ll}
0\quad&\mbox{if }p<p^{\rm aB}_{\rm c},\\
1\quad&\mbox{if and only if }p>p^{\rm aB}_{\rm c}.
\ea\right.
\ec
\et

\begin{figure}
\begin{center}
\inputtikz{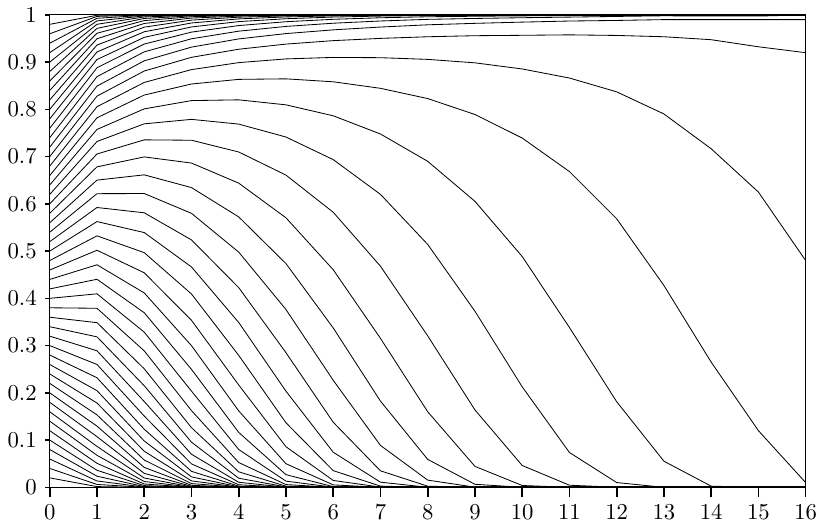}
\end{center}
\caption{The probability $P^{\rm Ab}_n(p)$ that Bob has a winning strategy in the game ${\rm Ab}_n(p)$, as a function of $n$, for different values of $p$. The data suggest a critical value $p_{\rm c}\approx 0.71$. They also suggest that $\lim_{n\to\infty}P^{\rm Ab}_n(p_{\rm c})$ is either one or very close to one.}\label{fig:treesim}
\end{figure}

Contrary to (\ref{Pearl}), we have not been able to determine the limit behaviour as $n\to\infty$ of $P^{\rm Ab}_n(p)$ at $p=p^{\rm Ab}_{\rm c}$, but we know that the limit cannot be zero. Numerical data suggest the limit is either one or close to one, see Figure~\ref{fig:treesim}. For $P^{\rm aB}_n(p)$, the situation is similar except that the roles of 0 and 1 have been interchanged. Numerical simulations suggest that $p^{\rm Ab}_{\rm c}\approx 0.71$ and $p^{\rm aB}_{\rm c}\approx 0.16$. We have the following rigorous bounds.

\bp[Bounds on the thresholds]
One\label{P:bounds} has $1/2\leq p^{\rm Ab}_{\rm c}\leq 7/8$ and $1/16\leq p^{\rm aB}_{\rm c}\leq\ha(3-\sqrt{5})$.
\ep

Recall that $\ha(3-\sqrt{5})=p^{\rm AB}_{\rm c}$ is the threshold from (\ref{Pearl}), so the bound $p^{\rm aB}_{\rm c}\leq\ha(3-\sqrt{5})$ follows from the inequality $P^{\rm aB}_n(p)\geq P^{\rm AB}_n(p)$ of Proposition~\ref{P:treecomp}.

The bound $1/2\leq p^{\rm Ab}_{\rm c}$ comes from the following observation. We can partition the possible outcomes of the game into $2^n$ sets of $n+1$ elements each, such that Alice's moves determine in which set the outcome lies, and then Bob's moves determine the precise outcome within each set. The probability that all outcomes in a given set are a win for Alice is $(1-p)^{n+1}$. Since there are $2^n$ sets, we see that if $p<1/2$, then in the limit $n\to\infty$, with high probability, there is at least one set in which all elements are a win for Alice. This means that Alice has a particularly simple winning strategy: by choosing her moves in advance, she can make sure she wins without even having to react to Bob's moves.

The proofs of the bounds $p^{\rm Ab}_{\rm c}\leq 7/8$ and $1/16\leq p^{\rm aB}_{\rm c}$ are a bit more complicated. We will use a Peierls argument first invented by Toom \cite{Too80} and further developed in \cite{SST22}. Toom's Peierls argument was invented to study perturbations of monotone cellular automata started in the all-one initial state. In Toom's original work, the time evolution is subject to random perturbations, but as was recently demonstrated by Hartarsky and Szab\'o in the context of bootstrap percolation \cite{HS22}, his method can also deal with perturbations of the initial state. To prove that the threshold in Theorem~\ref{T:threshold} is sharp, we will use an inequality from the theory of noise sensitivity first derived by Bourgain, Kahn, Kalai, Katznelson, and Linial \cite{BK+92}.

For the game ${\rm ab}_n(p)$ our results are much less complete than for the games ${\rm Ab}_n(p)$ and ${\rm aB}_n(p)$. Using Toom's Peierls argument, we are able to prove the following result, however.

\bp[Bounds for small and large parameters]
One\label{P:ab} has
\be
P^{\rm ab}_n(p)\asto{n}\left\{\ba{ll}
0\quad&\mbox{if }p<1/64,\\
1\quad&\mbox{if }p>15/16.
\ea\right.
\ee
\ep

The paper is organised as follows. In the remainder of the introduction, we discuss the graph structure of our games, we give an alternative interpretation of the quantities $P^{\rm AB}_n(p)$, $P^{\rm Ab}_n(p)$, $P^{\rm aB}_n(p)$, and $P^{\rm ab}_n(p)$ in terms of a cellular automaton, and we discuss related literature. Detailed proofs are deferred till Section~\ref{S:proof}.

\subsection{The graph structure of the games}\label{S:minmax}

To describe the graph structure of our games of interest, we first introduce some general notation. Recall that a directed graph is a pair $(G,\vec E)$ where $G$ is a set whose elements are called vertices and $\vec E$ is a subset of $G\times G$ whose elements are called directed edges. By a slight abuse of notation, we will sometimes use $G$ as a shorthand for $(G,\vec E)$, when it is clear how the set of directed edges is defined. We say that $(G,\vec E)$ is finite if $G$ is a finite set. For $v,w\in G$ we write $v\leadsto w$ if there exist $v=v_0,\ldots,v_n=w$ such that $(v_{k-1},v_k)\in\vec E$ for all $1\leq k\leq n$. We let ${\rm dist}(v,w)$ denote the smallest integer $n\geq 0$ for which such $v_0,\ldots,v_n$ exist. We say that $G$ is \emph{acyclic} if there do not exist $v,w\in G$ with $v\neq w$ such that $v\leadsto w\leadsto v$. For $v\in G$, we adopt the notation $\Oi(v):=\{w\in V:(v,w)\in\vec E\}$ and we write
\be
\ring{G}:=\big\{v\in V:\Oi(v)\neq\emptyset\big\}\quand\pa G:=\big\{v\in V:\Oi(v)=\emptyset\big\}.
\ee
Quite generally, we can describe the graph structure of a finite turn-based game played by two players by means of a quadruple $(G,\vec E,0,\tau)$ where
\begin{enumerate}
\item $(G,\vec E)$ is a finite acyclic directed graph,
\item $0\in G$ satisfies $0\leadsto v$ for all $v\in G$,
\item $\tau:\ring{G}\to\{\Ali,\Bob\}$ is a function.
\end{enumerate}
For the sake of the exposition, let us call such a quadruple a \emph{game-graph}. Slightly abusing the notation, we will sometimes write $G$ as shorthand for $(G,\vec E,0,\tau)$, when it is clear how $G$ is equipped with the structure of a game-graph. We call $0$ the \emph{root}. We interpret $G$ as the set of possible states of the game, $0$ as the state at the beginning of the game, $\pa G$ as the set of possible outcomes of the game, $\ring{G}$ as the non-final states of the game, and $\tau$ as a function that tells us whose turn it is in each non-final state of the game. We write
\be
\ring{G}^\Ali:=\big\{v\in\ring{G}:\tau(v)=\Ali\big\}
\quand
\ring{G}^\Bob:=\big\{v\in\ring{G}:\tau(v)=\Bob\big\}.
\ee
A (pure) strategy for Alice is a function
\be
\ring{G}^\Ali\ni v\mapsto\sig_\Ali(v)\in\Oi(v),
\ee
which has the interpretation that in the state $v$, Alice always plays the move that makes the game progress to the state $\sig_\Ali(v)$. Likewise, a strategy for Bob is a function $\ring{G}^\Bob\ni v\mapsto\sig_\Bob(v)\in\Oi(v)$. We let $\Si_\Ali$ (resp.\ $\Si_\Bob$) denote the set of all strategies for Alice (resp.\ Bob). Given a strategy $\sig_\Ali\in\Si_\Ali$ of Alice and a strategy $\sig_\Bob\in\Si_\Bob$ of Bob, there is a unique sequence of vertices $v_0,\ldots,v_n\in G$ with $v_0=0$ and $v_n\in\pa G$ such that for each $0\leq k<n$ one has $v_{k+1}=\sig_\Ali(v_k)$ if $v_k\in\ring{G}^\Ali$ and $v_{k+1}=\sig_\Bob(v_k)$ if $v_k\in\ring{G}^\Bob$. Then
\be\label{out}
o(\sig_\Ali,\sig_\Bob):=v_n
\ee
is the outcome of the game if Alice plays strategy $\sig_\Ali$ and Bob plays strategy $\sig_\Bob$. Given a function $x:\pa G\to\{0,1\}$ that assigns a winner to each possible outcome, with $x(v)=0$ meaning that Alice is the winner and $x(v)=1$ meaning that Bob is the winner, we say that a strategy $\sig_\Ali\in\Si_\Ali$ is a \emph{winning strategy} for Alice given $x$ if
\be
x\big(o(\sig_\Ali,\sig_\Bob)\big)=0\qquad\forall\sig_\Bob\in\Si_\Bob.
\ee
Winning strategies for Bob are defined similarly (with $1$ instead of $0$). To find out who has a winning strategy given $x$, we can work our way back from the set of all possible outcomes to the initial state of the game. Given a function $x:\pa G\to\{0,1\}$ that assigns a winner to each possible outcome, we can define $\ov x:\pa G\cup\ring{G}\to\{0,1\}$ by setting $\ov x(v):=x(v)$ for $v\in\pa G$ and then defining inductively
\be\label{minmax}
\ov x(v):=\left\{\ba{ll}
\dis\bigwedge_{w\in\Oi(v)}\ov x(w)\quad&\mbox{if }v\in\ring{G}^\Ali,\\[5pt]
\dis\bigvee_{w\in\Oi(v)}\ov x(w)\quad&\mbox{if }v\in\ring{G}^\Bob.
\ea\right.
\ee
Then it is not hard to see that Alice (resp.\ Bob) has a winning strategy if and only if $\ov x(0)=0$ (resp.\ $\ov x(0)=1$). We let $L:\{0,1\}^{\pa G}\to\{0,1\}$ denote the monotone Boolean function defined as
\be\label{Fdef}
L(x):=\ov x(0)\quad\mbox{where $\ov x$ solves (\ref{minmax}) with $\ov x=x$ on $\pa G$.}
\ee
Then $L(x)=0$ (resp.\ $=1$) if Alice (resp.\ Bob) has a winning strategy given $x$.

The game-graphs of the games that we are interested in have a special structure: the players alternate turns and play a fixed number of turns. Game graphs with various numbers of turns are consistent, in the sense that the game-graph of the game with $n$ turns can be obtained by truncating a certain infinite graph at height $n$. Moreover, these infinite graphs have a sort of product structure, in the sense that they are the ``product'' of two directed graphs associated with the individual players. To describe this, it will be useful to introduce some more notation.

Let us define a \emph{decision graph} to be a triple $(D,\vec F,0)$ such that
\begin{enumerate}
\item $(D,\vec F)$ is a directed graph,
\item $0\in D$ satisfies $0\leadsto w$ for all $w\in D$,
\item $|w|=|v|+1$ for all $(v,w)\in\vec F$, with $|w|:={\rm dist}(0,w)$ $(w\in D)$,
\item $\Oi(w)\neq\emptyset$ for all $w\in D$.
\end{enumerate}
Note that (iii) implies that $(D,\vec F)$ is acyclic. For a decision graph $D$, we introduce the notation
\be\ba{c}\label{decgr}
\dis\pa_n D:=\big\{w\in D:|w|=n\big\},\\[5pt]
\dis\li D\re_n:=\big\{w\in D:|w|<n\big\},\quand
[D]_n:=\big\{w\in D:|w|\leq n\big\}.
\ec
Given a decision graph $D$ and integer $n\geq 1$, we can define a game-graph
\be\label{GnDD1}
G_n(D)=(G_n,\vec E_n,0,\tau)
\ee
by setting
\be\label{GnDD2a}
G_n:=[D]_n\quand
\vec E_n:=\big\{(v,w)\in G_n\times G_n:(v,w)\in\vec F\big\},
\ee
which has the effect that $\ring G_n=\li D\re_n$ and $\pa G_n=\pa_nD$, and then defining
\be\label{GnDD2b}
\tau(v):=\left\{\ba{ll}
\dis\Ali\quad&\dis\mbox{if }|v|\mbox{ is even,}\\[5pt]
\dis\Bob\quad&\dis\mbox{if }|v|\mbox{ is odd.}
\ea\right.\qquad(v\in\ring G_n).
\ee
This corresponds to a game where Alice and Bob take turns to play $n$ moves together, and Alice starts. For the games we are interested in, $D$ has a product structure, in the sense that elements of $D$ are of the form $(a,b)$ where $a$ and $b$ record all that is relevant of Alice's and Bob's moves up to this point in the game. To formalise this, given two decision graphs $(D^i,\vec F^i,0^i)$ $(i=1,2$), we define a third decision graph $D^1\ltimes D^2:=(D,\vec F,0)$ by setting $0:=(0^1,0^2)$,
\be\label{ltimes}\left.\bac
\dis\pa_{2n}D&:=&\dis\big\{(a,b):a\in\pa_nD^1,\ b\in\pa_nD^2\big\},\\[5pt]
\dis\pa_{2n+1}D&:=&\dis\big\{(a,b):a\in\pa_{n+1}D^1,\ b\in\pa_nD^2\big\},
\ea\right\}\quad(n\in\N),
\ee
and
\bc\label{ltimesF}
\dis\vec F&:=&\dis\bigcup_{n=0}^\infty\{\big((a,b),(a',b)\big):(a,b)\in\pa_{2n}D,\ (a,a')\in\vec F^1\big\}\\[5pt]
&&\dis\cup\bigcup_{n=0}^\infty\{\big((a,b),(a,b')\big):(a,b)\in\pa_{2n+1}D,\ (b,b')\in\vec F^2\big\}.
\ec
The game-graphs of the games ${\rm AB}_n(p)$, ${\rm Ab}_n(p)$, ${\rm aB}_n(p)$, and ${\rm ab}_n(p)$ from the previous subsection are of the special form $G_{2n}(D^1\ltimes D^2)$ for a suitable choice of $D^1$ and $D^2$. In fact, only two different choices for $D^1$ and $D^2$ occur, which we now describe.

\begin{figure}
\begin{center}
\includegraphics[scale=0.7]{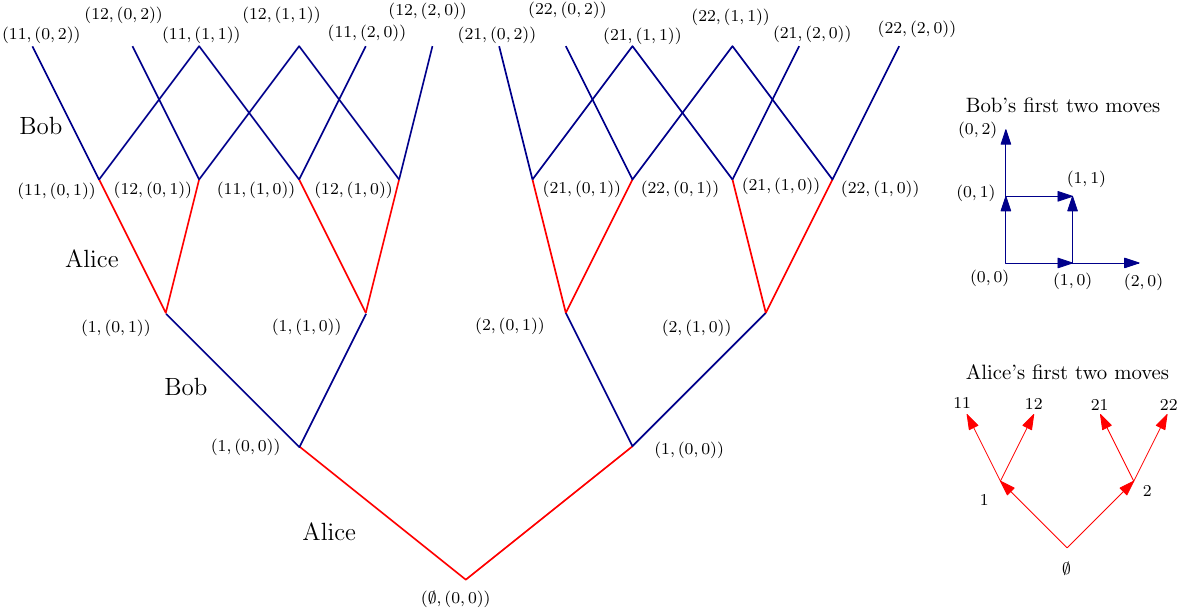}
\end{center}
\caption{The game-graph $G_4(\T\ltimes\N^2)$ of the game ${\rm Ab}_2(p)$.}\label{fig:treeBA}
\end{figure}

Let $\T$ denote the space of all finite words $\ibf=i_1\cdots i_l$ $(l\geq 0)$ made up from the alphabet $\{1,2\}$. We call $|\ibf|:=l$ the \emph{length} of $\ibf$ and we let $\wurz$ denote the word of length $l=0$. If $\ibf=i_1\cdots i_l$ and $\jbf=j_1\cdots j_m$, then we let $\ibf\jbf=i_1\cdots i_lj_1\cdots j_m$ denote the concatenation of $\ibf$ and $\jbf$. We view $\T$ as a decision graph $(D,\vec F,0)$ with vertex set $D:=\T$, edge set
\be
\vec F:=\big\{(\ibf,\ibf j):\ibf\in\T,\ j\in\{1,2\}\big\},
\ee
and root $0:=\wurz$. Note that $|\ibf|={\rm dist}(\wurz,\ibf)$.

Let $\N^2$ denote the set of all pairs $(i,j)$ of nonnegative integers. We view $\N^2$ as a decision graph $(D,\vec F,0)$ with vertex set $D:=\N^2$, edge set
\be
\vec F:=\big\{\big((i,j),(i+1,j)\big),\big((i,j),(i,j+1)\big):(i,j)\in\N^2\big\},
\ee
and root $0:=(0,0)$. Note that $\big|(i,j)\big|:=i+j={\rm dist}\big((0,0),(i,j)\big)$.

With this notation, the game-graphs of the games ${\rm AB}_n(p)$, ${\rm Ab}_n(p)$, ${\rm aB}_n(p)$, and ${\rm ab}_n(p)$ are given by
\be\label{gagr}
G_{2n}(\T\ltimes\T),\quad G_{2n}(\T\ltimes\N^2),\quad G_{2n}(\N^2\ltimes\T),\quand G_{2n}(\N^2\ltimes\N^2),
\ee
respectively. We let
\be\label{FAB}
L^{\rm AB}_{2n},\quad L^{\rm Ab}_{2n},\quad L^{\rm aB}_{2n},\quand L^{\rm ab}_{2n}
\ee
denote the corresponding monotone Boolean functions, defined as in (\ref{Fdef}). Then
\be\label{Pminmax}
P^{\rm AB}_n(p)=\P\big[L^{\rm AB}_{2n}(X^p)=1\big]\quad\mbox{where}\quad
\big(X^p(v)\big)_{v\in\pa G_{2n}}\mbox{ are i.i.d.\ with intensity }p,
\ee
and the same holds with ${\rm AB}$ replaced by ${\rm Ab}$, ${\rm aB}$, or ${\rm ab}$. Here $\pa G_{2n}(\T\ltimes\T)=\pa_{2n}(\T\ltimes\T)=(\pa_n\T)\times(\pa_n\T)$ and similarly in the other three cases.

\begin{figure}
\begin{center}
\inputtikz{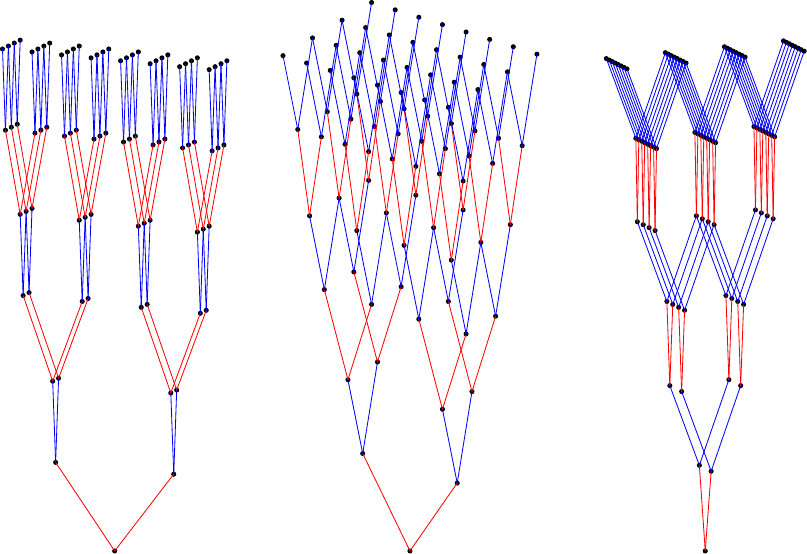}
\end{center}
\caption{Three projections of the game-graph $G_6(\T\ltimes\N^2)$ of the game ${\rm Ab}_3(p)$. Seen from one side (picture on the left), the graph looks like the binary tree of height three $[\T]_3$ (in red), with vertical edges added (in blue). Seen from the other side (picture on the right), the graph looks like the game-graph $[\N^2]_3$ (in blue), with vertical edges added (in red). In Figures \ref{fig:Toom} and \ref{fig:count} below, we make use of such projections of $G_6(\T\ltimes\N^2)$ and $G_4(\N^2\ltimes\N^2)$.}\label{fig:gamegraph}
\end{figure}

We observe that $\T\ltimes\T\cong\T$, i.e., the decision graph $\T\ltimes\T$ has the structure of a binary tree. This means that in the game ${\rm AB}_n(p)$, to each possible outcome leads a unique game history. On the other hand, in the games ${\rm Ab}_n(p)$, ${\rm aB}_n(p)$, and ${\rm ab}_n(p)$, different game histories can lead to the same outcome. See Figures \ref{fig:treeBA} and \ref{fig:gamegraph} for pictures of the game-graphs of the games ${\rm Ab}_2(p)$ and ${\rm Ab}_3(p)$.

\subsection{Cellular automata}\label{S:cell}

In Subsection~\ref{S:intro} we gave a somewhat informal description of the functions $P^{\rm AB}_n$, $P^{\rm Ab}_n$, $P^{\rm aB}_n$, and $P^{\rm ab}_n$ in terms of winning strategies. In formula (\ref{Pminmax}) of Subsection~\ref{S:minmax} we gave an alternative description of these functions in terms of Boolean functions with i.i.d.\ input. In the present subsection we will give yet another representation of these functions in terms of cellular automata.

\begin{figure}
\begin{center}
\inputtikz{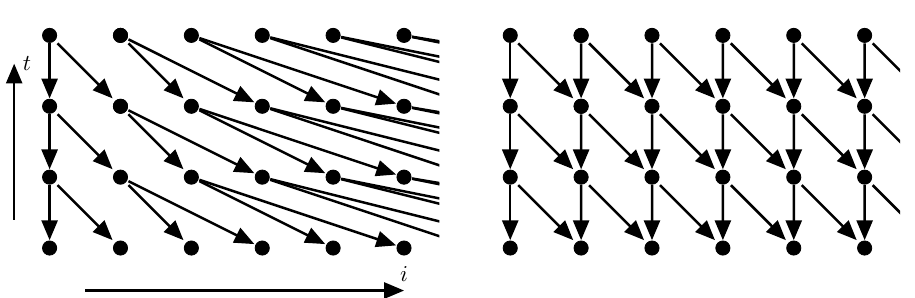}
\end{center}
\caption{Time dependencies for rules ${\rm A}$ and ${\rm a}$. In the picture on the left, in line with rule ${\rm A}$, we have drawn arrows from each point $(i,t)\in\N\times\N$ to the points $(2i,t-1)$ and $(2i+1,t-1)$, which has the result that the set of points that can be reached from $(0,t)$ has the structure of a binary tree of depth $t$. In the picture on the right, in line with rule ${\rm a}$, we have drawn arrows from each point $(i,t)\in\N\times\Z$ to the points $(i,t-1)$ and $(i+1,t-1)$, which has the result that the set of points that can be reached from $(0,t)$ has the structure of the game-graph $G_t(\N^2)$. In a similar way, one can equip $\N^2\times\N$ with the structure of a directed graph that corresponds to applying at even times one of the rules ${\rm A}$ or ${\rm a}$ in one spatial direction, and at odd times one of the rules ${\rm B}$ or ${\rm b}$ in the other spatial direction. In this case, the ``genealogy'' of a space-time point $(0,0,t)$ is described by one of the game-graphs $G_t(\T\ltimes\T)$, $G_t(\T\ltimes\N^2)$, $G_t(\N^2\ltimes\T)$, or $G_t(\N^2\ltimes\N^2)$.}\label{fig:cellaut}
\end{figure}

Let $\{0,1\}^{\N^2}$ be the set of functions $x:\N^2\to\{0,1\}$. Given a random variable $X_0$ with values in $\{0,1\}^{\N^2}$, we will be interested in a stochastic process $(X_t)_{t\in\N}$ that evolves in a deterministic way. At even times $t>0$, we define $X_t$ in terms of $X_{t-1}$ according to one of the following two rules
\be\ba{r@{\quad}l}\label{Arules}
{\rm A}.&\dis X_t(i,j)=X_{t-1}(2i,j)\wedge X_{t-1}(2i+1,j),\\[5pt]
{\rm a.}&\dis X_t(i,j)=X_{t-1}(i,j)\wedge X_{t-1}(i+1,j),
\ea\qquad\big((i,j)\in\N^2\big),
\ee
and at odd times $t>0$, we define $X_t$ according to one of the rules
\be\ba{r@{\quad}l}\label{Brules}
{\rm B}.&\dis X_t(i,j)=X_{t-1}(i,2j)\vee X_{t-1}(i,2j+1),\\[5pt]
{\rm b}.&\dis X_t(i,j)=X_{t-1}(i,j)\vee X_{t-1}(i,j+1),
\ea\qquad\big((i,j)\in\N^2\big).
\ee
The four possible combinations yield four possible ways to define a dynamic, which we denote by ${\rm AB}$, ${\rm Ab}$, ${\rm aB}$, and ${\rm ab}$. Let $X^p_0=\big(X^p_0(i,j)\big)_{(i,j)\in\N^2}$ be i.i.d.\ $\{0,1\}$-valued random variables with $\P[X^p_0(i,j)=1]=p$ $(i,j\in\N)$, and let $(X^{{\rm AB},p}_t)_{t\in\N}$ be defined according to the rules ${\rm A}$ and ${\rm B}$ with initial state $X^{{\rm AB},p}_0:=X^p_0$. We claim that
\be\label{Pcell}
P^{\rm xx}_n(p)=\P\big[X^{{\rm xx},p}_{2n}(0,0)=1\big]\quad \mbox{where }{\rm xx}={\rm AB},\ {\rm Ab},\ {\rm aB},\mbox{ or }{\rm ab}.
\ee
We only sketch the proof for the combination of rules ${\rm Ab}$. The main idea is to look at the ``genealogy'' of a space-time point $(i,j,2n)$. By rule ${\rm A}$, we have
\be
X_{2n}(i,j)=X_{2n-1}(2i,j)\wedge X_{2n-1}(2i+1,j),
\ee
where by rule ${\rm b}$,
\bc
\dis X_{2n-1}(2i,j)&=&\dis X_{2n-2}(2i,j)\vee X_{2n-2}(2i,j+1),\\[5pt]
\dis X_{2n-1}(2i+1,j)&=&\dis X_{2n-2}(2i+1,j)\vee X_{2n-2}(2i+1,j+1).
\ec
Continuing in this way, one can work out how $X_{2n}(0,0)$ is defined in terms of $\big(X_0(i,j)\big)_{(i,j)\in\N^2}$. A careful inspection of the rules (see Figure~\ref{fig:cellaut}) shows that this leads precisely to the structure of the game-graph $G_{2n}(\T\ltimes\N^2)$, so that (\ref{Pcell}) follows from (\ref{Pminmax}). The proofs for the other combinations of rules are similar.

Because of the tree structure of the decision graph $\T\ltimes\T$, it is easy to see that for the cellular automaton $\big(X^{{\rm AB},p}_t\big)_{t\in\N}$, one has that
\be
\big(X^{{\rm AB},p}_t(i,j)\big)_{(i,j)\in\N^2}\quad\mbox{are i.i.d.\ for all }t\in\N.
\ee
Because of this, the analysis of the game ${\rm AB}_n(p)$ is relatively simple. For the cellular automaton $\big(X^{{\rm Ab},p}_t\big)_{t\in\N}$, we will prove in Lemma~\ref{L:column} below that it is still true that
\be\ba{l}\label{column}
\big(X^{{\rm Ab},p}_t(i,\cdot)\big)_{i\in\N}\quad\mbox{are i.i.d.\ for all }t\in\N\\[5pt]
\quad\mbox{where}\quad X^{{\rm Ab},p}_t(i,\cdot):=\big(X^{{\rm Ab},p}_t(i,j)\big)_{j\in\N}\qquad(i\in\N),
\ec
but in the $j$-direction dependencies develop. Because of this, the analysis of the game ${\rm Ab}_n(p)$ is more difficult and it seems unlikely that an explicit formula for the critial value $p^{\rm Ab}_{\rm c}$ can be found. The situation for the game ${\rm aB}_n(p)$ is similar, while for the cellular automaton $\big(X^{{\rm ab},p}_t\big)_{t\in\N}$ independence is lost in both directions.

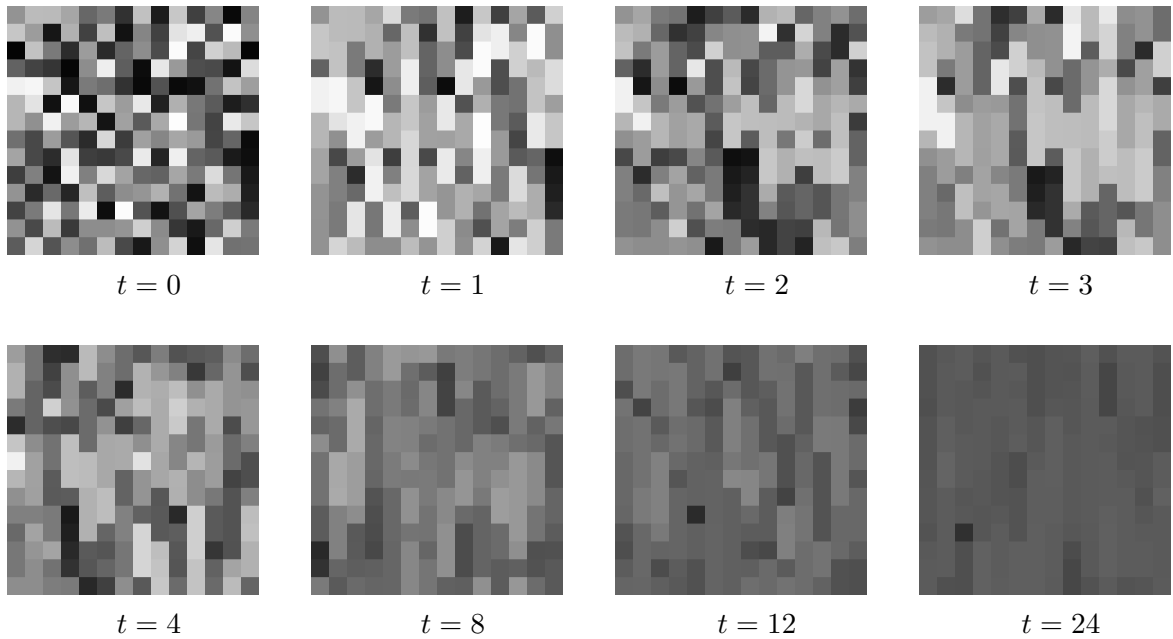
\begin{figure}
\begin{center}
\begin{tikzpicture}
\node[above right] at (0,0) {\includegraphics[width=3.3cm]{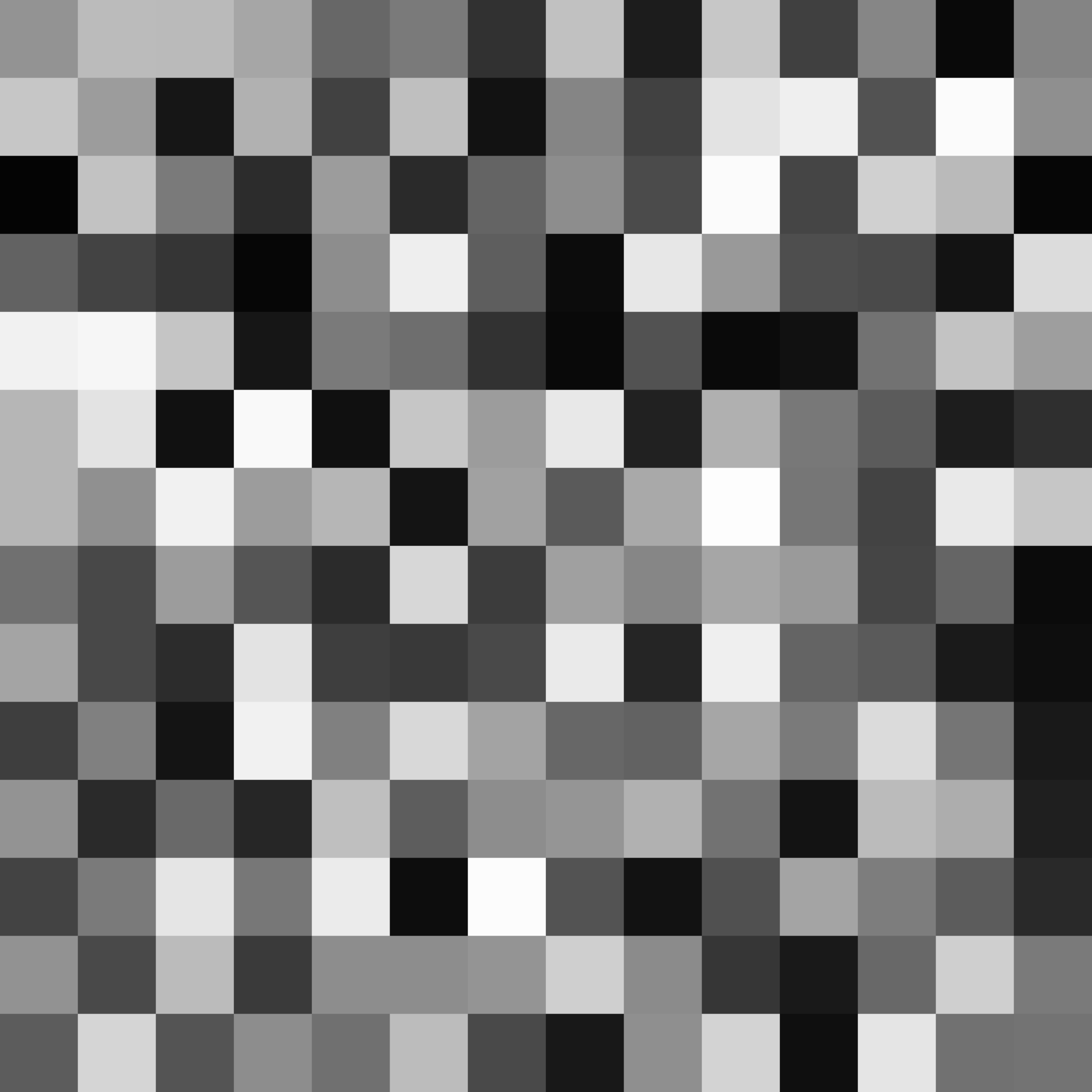}};
\node[above right] at (4,0) {\includegraphics[width=3.3cm]{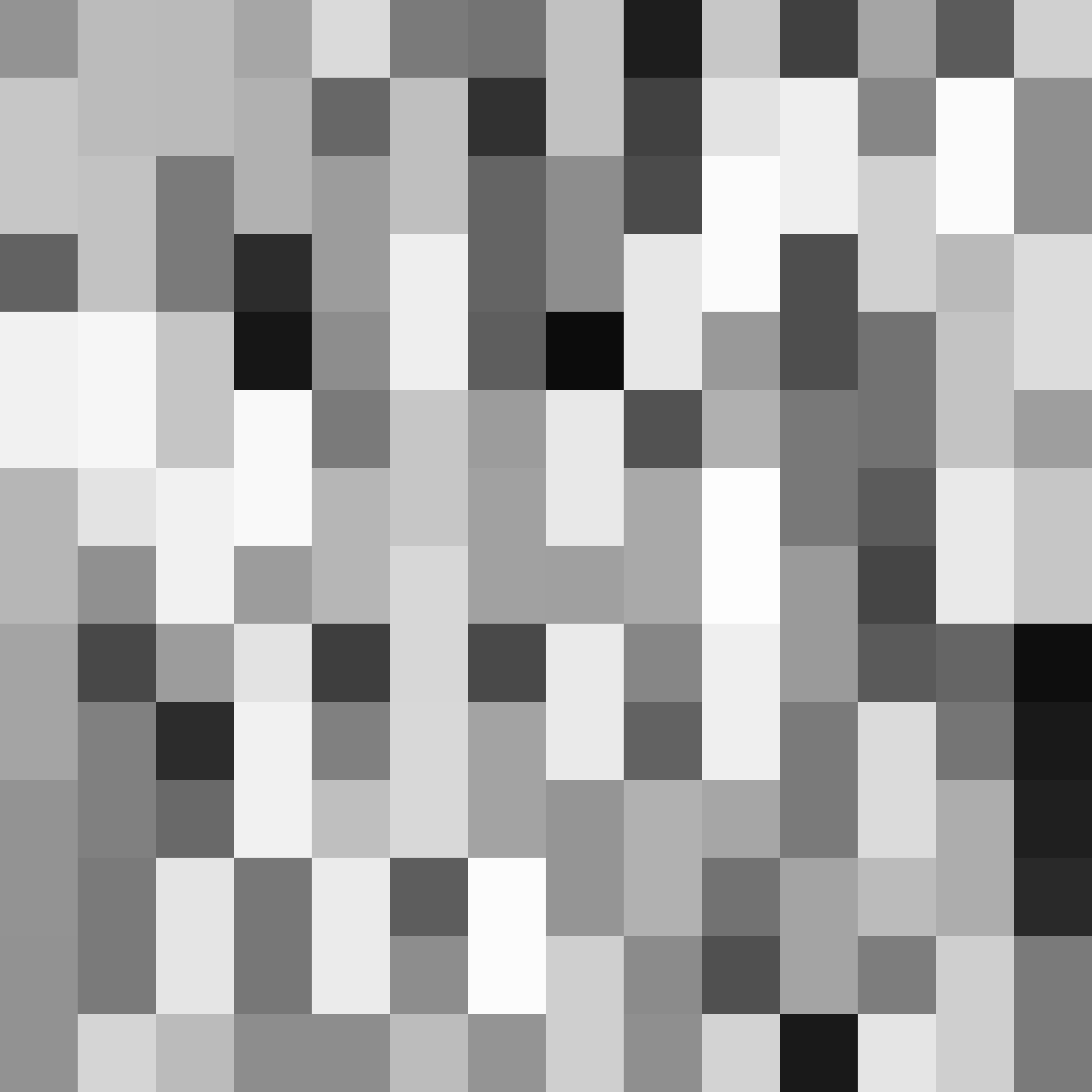}};
\node[above right] at (8,0) {\includegraphics[width=3.3cm]{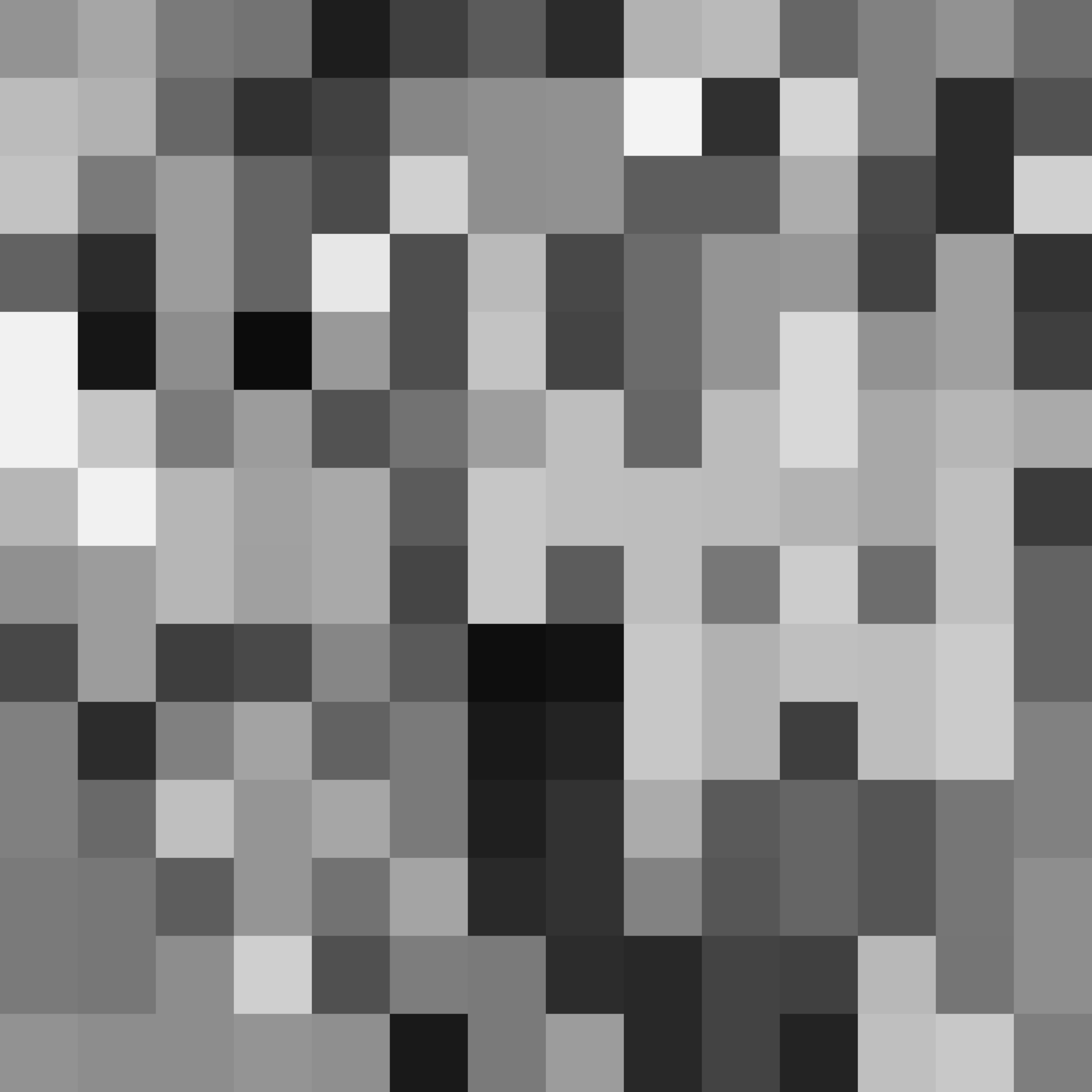}};
\node[above right] at (12,0) {\includegraphics[width=3.3cm]{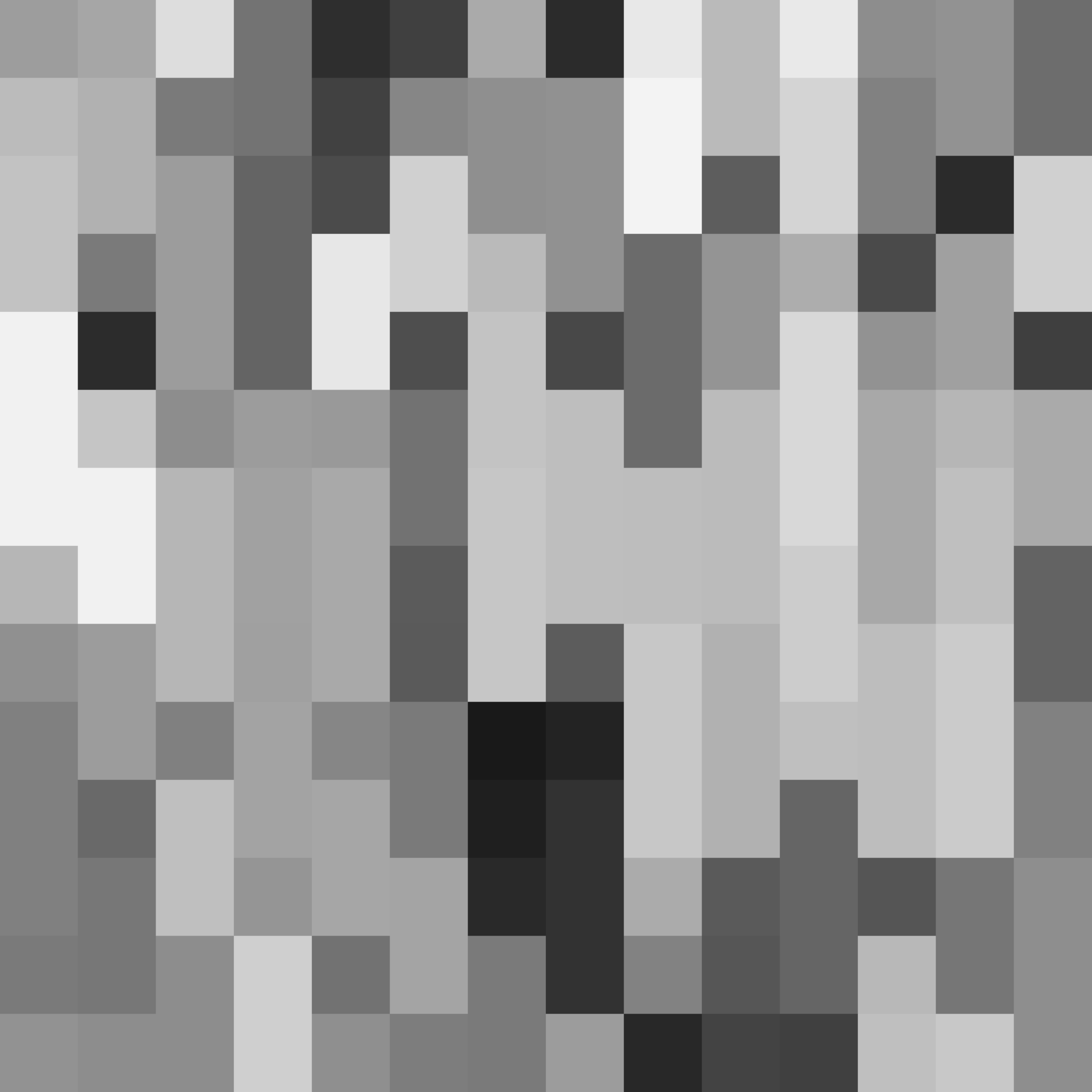}};
\node[below] at (2,0) {$t=0$};
\node[below] at (6,0) {$t=1$};
\node[below] at (10,0) {$t=2$};
\node[below] at (14,0) {$t=3$};
\node[above right] at (0,-4.5) {\includegraphics[width=3.3cm]{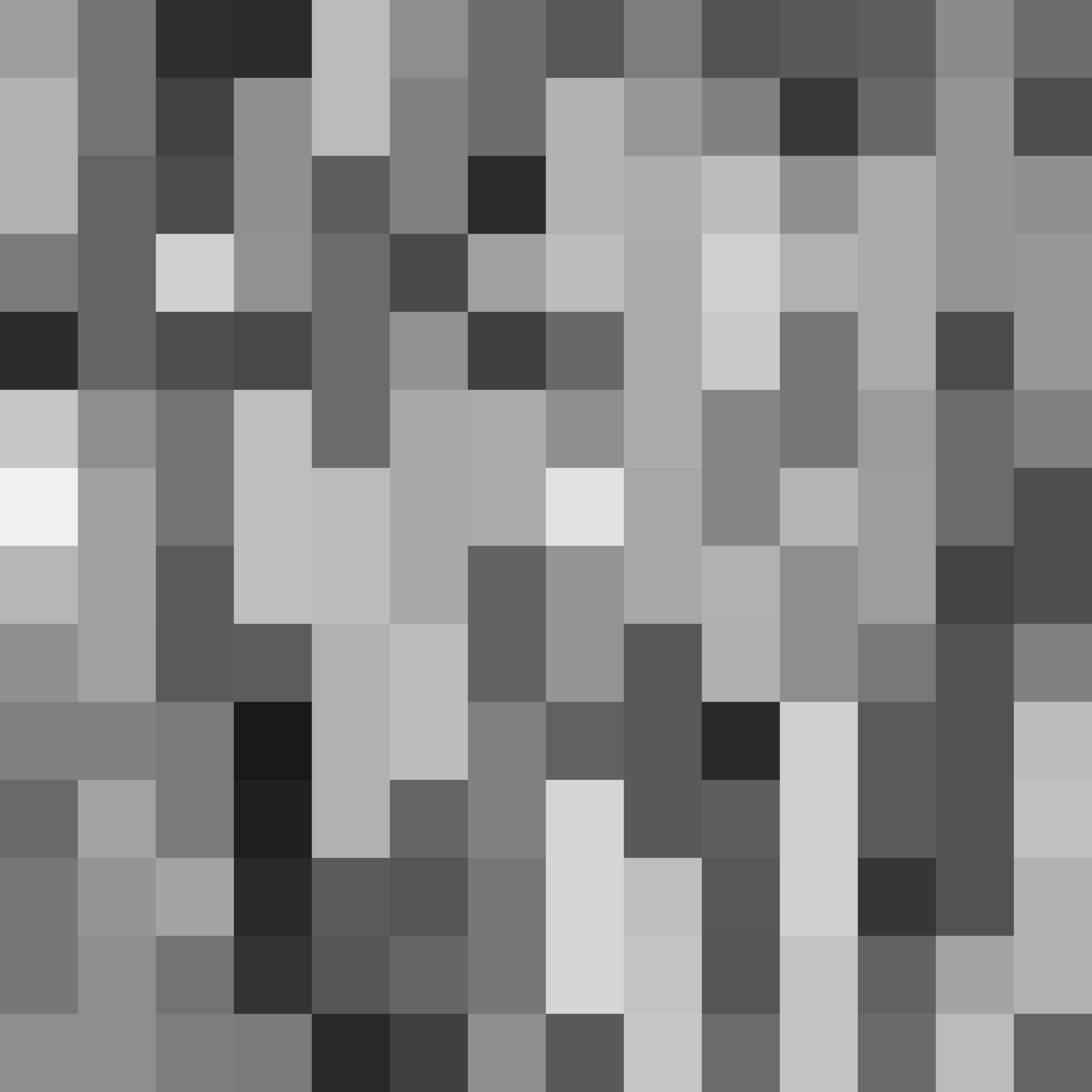}};
\node[above right] at (4,-4.5) {\includegraphics[width=3.3cm]{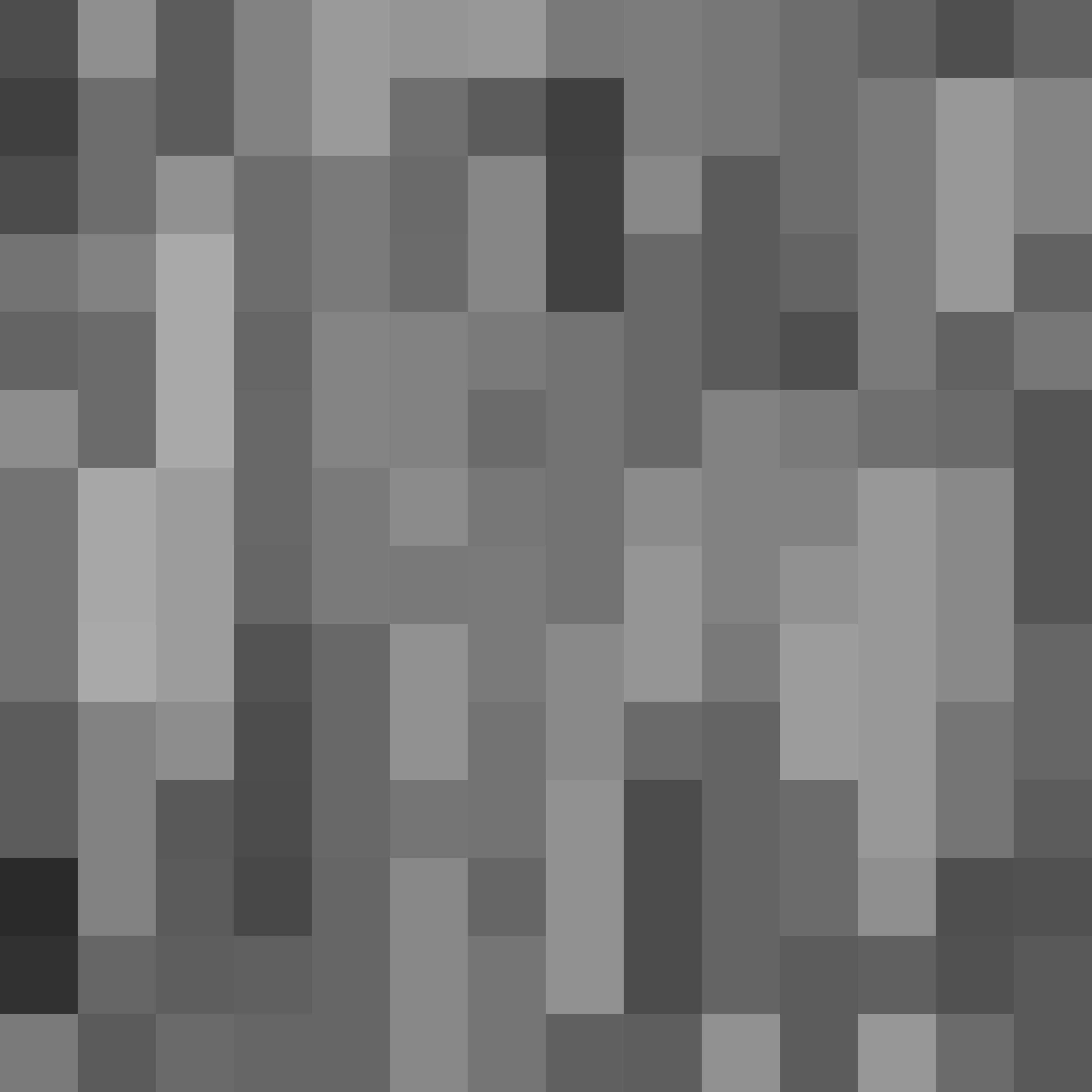}};
\node[above right] at (8,-4.5) {\includegraphics[width=3.3cm]{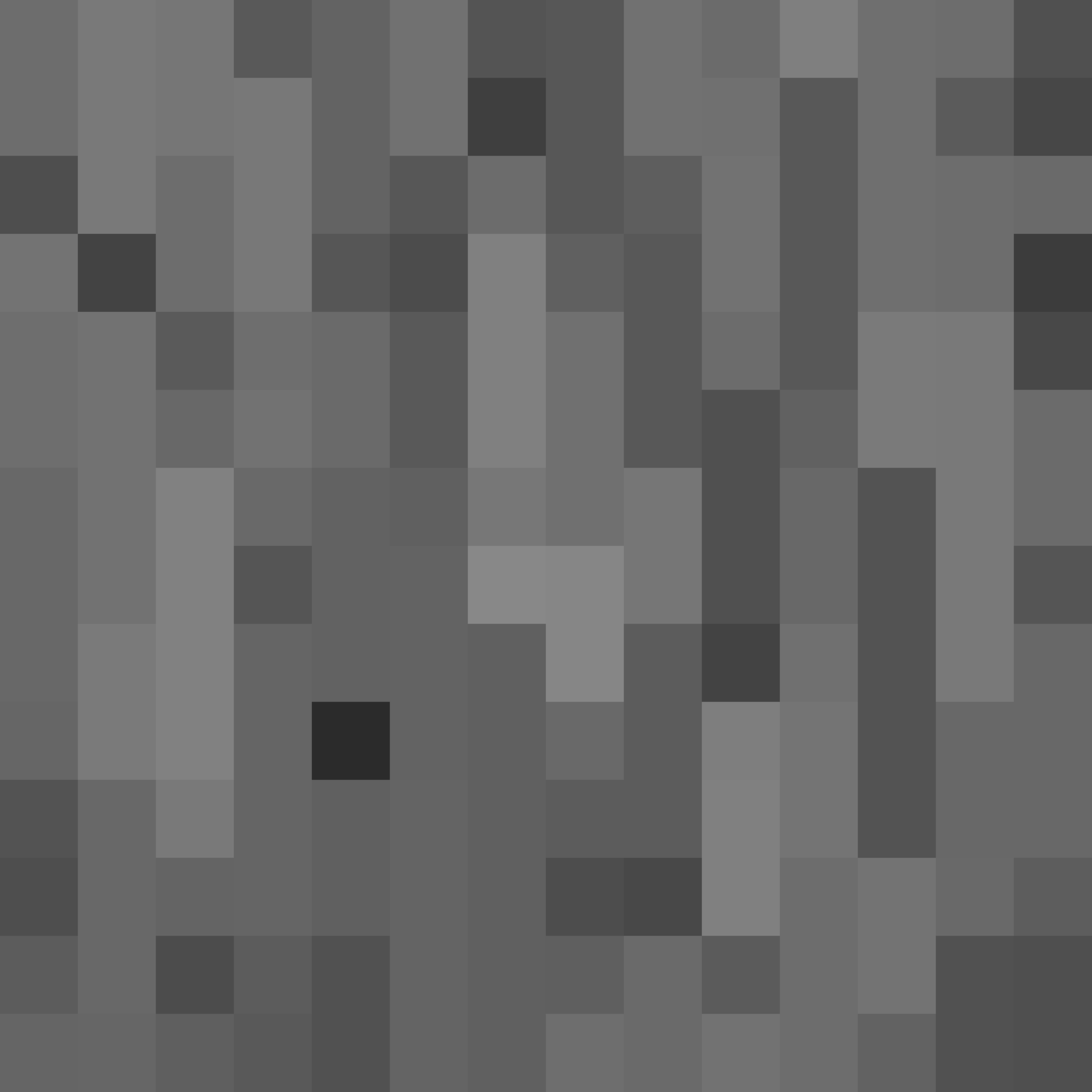}};
\node[above right] at (12,-4.5) {\includegraphics[width=3.3cm]{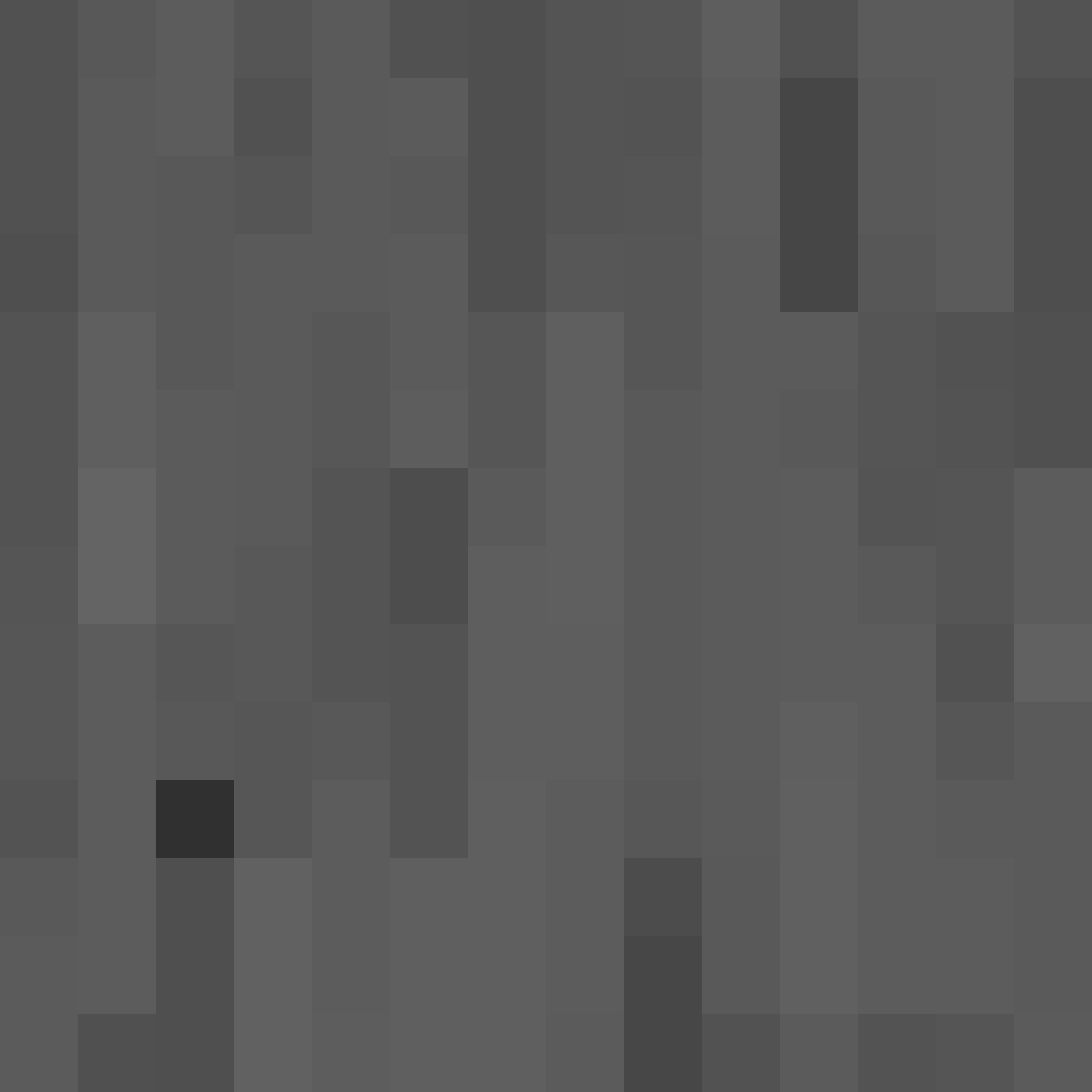}};
\node[below] at (2,-4.5) {$t=4$};
\node[below] at (6,-4.5) {$t=8$};
\node[below] at (10,-4.5) {$t=12$};
\node[below] at (14,-4.5) {$t=24$};
\end{tikzpicture}
\end{center}
\caption{Snapshots of the cellular automaton $(U^{\rm Ab}_t)_{t\in\N}$. Grayscales indicate a value between zero (white) and one (black). Initially, all lattice points are independent and uniformly distributed. At odd times, the value of $(i,j)$ is replaced by the maximum of the values at $(i,j)$ and $(i,j+1)$, and at even times, the value at $(i,j)$ is replaced by the minimum of the values at $(2i,j)$ and $(2i+1,j)$. Columns remain independent of each other but dependencies develop in the vertical direction. Our results imply that as time tends to infinity, the value of each lattice point tends to $p^{\rm Ab}_{\rm c}$ in probability.}
 \label{fig:Ucell}
\end{figure}

There exists a nice way of coupling cellular automata with different initial intensities $p$. Let $\big(U_0(i,j)\big)_{(i,j)\in\N^2}$ be i.i.d.\ uniformly distributed $[0,1]$-valued random variables. Then we can inductively define cellular automata with values in $[0,1]^{\N^2}$ by applying at even times $t>0$ one of the rules
\be\ba{r@{\quad}l}\label{AUrules}
{\rm A}.&\dis U_t(i,j)=U_{t-1}(2i,j)\vee U_{t-1}(2i+1,j),\\[5pt]
{\rm a.}&\dis U_t(i,j)=U_{t-1}(i,j)\vee U_{t-1}(i+1,j),
\ea\qquad\big((i,j)\in\N^2\big),
\ee
and at odd times $t>0$ one of the rules
\be\ba{r@{\quad}l}\label{BUrules}
{\rm B}.&\dis U_t(i,j)=U_{t-1}(i,2j)\wedge U_{t-1}(i,2j+1),\\[5pt]
{\rm b}.&\dis U_t(i,j)=U_{t-1}(i,j)\wedge U_{t-1}(i,j+1),
\ea\qquad\big((i,j)\in\N^2\big).
\ee
Note that compared to (\ref{Arules}) and (\ref{Brules}), the minimum and maximum operations have been interchanged. Defining $(U^{\rm AB}_t)_{t\in\N}$, $(U^{\rm Ab}_t)_{t\in\N}$, and so on in this way, it is straightforward to check that for each $p\in[0,1]$, setting
\be
X^{{\rm xx},p}_t(i,j):=\left\{\ba{ll}
1\quad&\mbox{if }U^{\rm xx}_t(i,j)\leq p,\\[5pt]
0\quad&\mbox{if }U^{\rm xx}_t(i,j)>p,
\ea\right.
\quad\mbox{where }{\rm xx}={\rm AB},\ {\rm Ab},\ {\rm aB},\mbox{ or }{\rm ab}.
\ee
defines a cellular automaton with values in $\{0,1\}^{\N^2}$ of the type described earlier. It follows therefore from (\ref{Pcell}) that
\be\label{PU}
P^{\rm xx}_n(p)=\P\big[U^{{\rm xx}}_{2n}(0,0)\leq p\big]\quad\mbox{where }{\rm xx}={\rm AB},\ {\rm Ab},\ {\rm aB},\mbox{ or }{\rm ab}.
\ee
See Figure~\ref{fig:Ucell} for a numerical simulation of the cellular automaton $(U^{\rm Ab}_t)_{t\in\N}$.

This coupling has a nice interpretation in terms of a game. Let $(G,\vec E,0,\tau)$ be a game-graph as defined in Subsection~\ref{S:minmax} and let $\big(U(v)\big)_{v\in\pa G}$ be i.i.d.\ uniformly distributed $[0,1]$-valued random variables, that have the interpretation that if the outcome of the game is $v$, then Alice receives the pay-out $U(v)$ and Bob receives the pay-out $1-U(v)$. If both players play optimally, then the game will end in some a.s.\ unique outcome $\vbf\in\pa G$, which is characterised by the fact that
\be
U(\vbf)=\sup_{\sig_\Ali\in\Si_\Ali}\inf_{\sig_\Bob\in\Si_\Bob}U\big(o(\sig_\Ali,\sig_\Bob)\big),
\ee
where as in (\ref{out}), $o(\sig_\Ali,\sig_\Bob)$ denotes the outcome of the game if Alice plays strategy $\sig_\Ali$ and Bob plays strategy $\sig_\Bob$. Specialising to the game-graphs from (\ref{gagr}), we can in this way define random games ${\rm AB}_n$, ${\rm Ab}_n$, ${\rm aB}_n$, and ${\rm ab}_n$ where the pay-out of Alice and Bob for each outcome is determined by i.i.d.\ uniformly distributed $[0,1]$-valued random variables. Then the random variable $1-U^{{\rm AB}}_{2n}(0,0)$ from (\ref{PU}) can be interpreted as Bob's pay-out in the game ${\rm AB}_n$ if both players play optimally, and similarly for the games ${\rm Ab}_n$, ${\rm aB}_n$, and ${\rm ab}_n$.

\subsection{Discussion and open problems}

The game-graph of Pearl's \cite{Pea80} original game ${\rm AB}_n(p)$ is the deterministic tree $G_{2n}(\T\ltimes\T)$. Several authors have studied variants of the game ${\rm AB}_n(p)$ where this tree is replaced by a random tree. Martin and Stasi\'nski in \cite{MS18} studied min-max games on a Galton-Watson tree, truncated at some height $2n$. Pemantle and Ward in \cite{PW06} studied a model on a binary tree where i.i.d.\ random variables attached to the nodes determine whether the minimum or maximum function should be applied. We refer to Broutin and Mailler \cite{BM18} for an overview of the literature on this type of problems.

For games with a large number of moves, it may be computationally unfeasible to find a winning strategy. This motivates the study of efficient algorithms that can play games near optimally. Indeed, this was the main motivation for Pearl \cite{Pea80} to study the game ${\rm AB}_n(p)$. There is an extensive literature on this topic. We cite \cite{DG95} for a somewhat older overview and \cite{SN15,MS18} for more recent contributions.

When trying to play a game near optimally, a natural approach is to find all possible states that can be reached after a fixed number $N$ of moves, then attach a utility to each possible state using some statistical procedure, and finally use minimaxing in the spirit of (\ref{minmax}) to compute utility values for all states that lie less than $N$ moves in the future. Surprisingly, it was discovered by Nau \cite{Nau80} that when the game-graph is a tree, such algorithms can sometimes behave pathologically, in the sense that increasing $N$ leads to worse, rather than better game playing. This sort of pathology has been the subject of much further research \cite{LBG12}.

An early suggestion, that has not been pursued much, is that this sort of pathologies may be resolved if different game histories can lead to the same outcome \cite{Nau83}. When the game no longer has the structure of a tree, however, independence is lost. For this reason, such games are harder to analyse and much less is known. There seems to be a general agreement, however, that introducing dependence may resolve game-playing pathology. For this reason, it would seem interesting to study game-playing algorithms for the games ${\rm Ab}_n(p)$ and ${\rm aB}_n(p)$ introduced in the present paper.

Our own motivation for studying these games comes from a different direction, which is related to the cellular automata from Subsection~\ref{S:cell}. Attractive spin systems have two special invariant laws, called the lower and upper invariant laws, that are the long-time limit laws starting from the all zero and all one initial states \cite[Thm~III.2.3]{Lig85}. If the dynamics of the spin system are invariant under automorphisms of the lattice (such as translations on $\Z^d$), then it is natural to ask if all invariant laws that are moreover invariant under automorphisms of the lattice are convex combinations of the lower and upper invariant laws. Frequently, one finds that mean-field models have an additional ``intermediate'' invariant law that lies between the lower and upper invariant laws, but very little is known about the existence of such intermediate invariant laws in a truly spatial setting. A notable exception are stochastic Ising models, where such intermediate invariant laws are known to exist if the lattice is a Cayley tree \cite[Section~12.2]{Geo11} but not on $\Z^d$. (The latter follows from \cite{Bod06} (in dimensions $d\geq 3$) and \cite{Hig81} (in dimension $d=2$) combined with \cite[Thm~IV.5.12]{Lig85}; see also the discussion on \cite[page~166]{FV18}).

The cellular automaton $X^{\rm AB}$ from Subsection~\ref{S:cell} is a mean-field model.\footnote{The cellular automaton $X^{\rm AB}$, started from an i.i.d. law, remains i.i.d. for all time. This is characteristic of mean-field models (sometimes called ``propagation of chaos'') and is closely related to the tree structure of the ``genealogy'' of $X^{\rm AB}$. Compare the relation between recursive tree processes and mean-field limits in \cite{MSS20}.} It has an intermediate invariant law, which is the product measure with intensity $p_{\rm c}^{{\rm AB}}$. Our original plan was to prove that the cellular automata $X^{\rm Ab}$ and $X^{\rm aB}$ also have an intermediate invariant law, but this is still unresolved. This is closely related to the question about the limit behaviour of $P^{\rm Ab}_n(p)$ at $p=p^{\rm Ab}_{\rm c}$. We have shown that this quantity does not tend to zero. It seems there are two plausible scenarios: either this quantity tends to one and the cellular automaton $X^{\rm Ab}$ has no intermediate invariant law, or it tends to a value strictly between zero and one and such an intermediate invariant law exists.

Another open problem concerns the size of the ``critical window'' where the function $P^{\rm Ab}_n$ changes from a value close to zero to a value close to one. In Theorem~\ref{T:sharp} below, we prove an upper bound of order $1/n$, but it is doubtful this is sharp, given that for the game ${\rm AB}_n(p)$ the size of the critical window decreases exponentially in $n$. For the mean-field game ${\rm AB}_n(p)$, a detailed analysis of the critical behaviour, which involves a nontrivial limit law, has been carried out in \cite{ADN05}.

There are lots of open problems concerning the game ${\rm ab}_n(p)$. Below Proposition~\ref{P:treecomp} we already mentioned the conjecture that $P^{ab}_n(p)\leq P^{aB}_n(p)$ and $P^{ab}_n(p)\geq P^{Ab}_n(p)$ for all $p\in[0,1]$ and $n\geq 1$. Another obvious conjecture is the existence of a critical value $p^{\rm ab}_{\rm c}$ such that
\be\label{pabc}
P^{\rm ab}_n(p)\asto{n}\left\{\ba{ll}
0\quad&\mbox{if }p<p^{\rm ab}_{\rm c},\\
1\quad&\mbox{if }p>p^{\rm ab}_{\rm c}.
\ea\right.
\ee
Our proof of Theorem~\ref{T:threshold} uses the special structure of the game-graphs associated with the games ${\rm Ab}_n(p)$ and ${\rm aB}_n(p)$ in two places: namely, in the preparatory Lemmas \ref{L:inflow} and  \ref{L:nul}  we use that the game-graphs have a product structure with one component a binary tree. Based on numerical simulations, we conjecture that a critical value $p^{\rm ab}_{\rm c}$ as in (\ref{pabc}) exists and that the cellular automaton $X^{\rm ab}$ started in product law with this intensity converges as time tends to infinity to a convex combination of the delta-measures on the all-zero and all-one states. Based on this, we conjecture that the cellular automaton $X^{\rm ab}$ has no intermediate invariant law.

\section{Proofs}\label{S:proof}

\subsection{An inequality for monotone Boolean functions}\label{S:bolineq}

In this subsection we prepare for the proof of Proposition~\ref{P:treecomp} by establishing a general comparison result for monotone Boolean functions applied to i.i.d.\ input.

A \emph{Boolean variable} is a variable that takes values in the set $\{0,1\}$. A \emph{Boolean function} is any function that maps a set of the form $\{0,1\}^S$, where $S$ is a finite set, into $\{0,1\}$. A Boolean function is \emph{monotone} if $x\leq y$ (pointwise) implies $L(x)\leq L(y)$. Let $L:\{0,1\}^S\to\{0,1\}$ be a monotone Boolean function. By definition, a \emph{one-set} of $L$ is a set $A\sub S$ such that $L(1_A)=1$, where $1_A$ denotes the indicator function of $A$, and a \emph{zero-set} of $L$ is a set $Z\sub S$ such that $L(1-1_Z)=0$. A one-set (respectively zero-set) is \emph{minimal} if it does not contain other one-sets (respectively zero-sets) as a proper subset. We let $\Ai^\up(L)$ and $\Zi^\up(L)$ denote the set of one-sets and zero-sets of $L$, respectively, and we denote the sets of minimal one-sets and zero-sets by $\Ai(L)$ and $\Zi(L)$, respectively.\footnote{This notation is motivated by the fact that $\Ai^\up(L)=\{A\sub S:\exists A'\in\Ai(L)\mbox{ s.t.\ }A'\sub A\}$ by the monotonicity of $L$ and also by the desire to have a simple notation for $\Ai(L)$ which is needed more often than $\Ai^\up(L)$.} We will prove the following result.

\bp[An inequality]
Let\label{P:compar} $S$ and $T$ be finite sets and let $L:\{0,1\}^S\to\{0,1\}$ be a monotone Boolean function. Let $\psi:S\to T$ be a function such that
\be\label{minint}
\big|\psi(A)\big|=|A|\qquad\forall A\in\Ai(L),
\ee
and let $X=(X(i))_{i\in S}$ and $Y=(Y(j))_{j\in T}$ by i.i.d.\ collections of Boolean random variables with intensity $p\in[0,1]$. Then
\be\label{compar}
\P\big[L(Y\circ\psi)=1\big]\leq\P\big[L(X)=1\big].
\ee
Similarly, if (\ref{minint}) holds with minimal one-sets replaced by minimal zero-sets, then (\ref{compar}) holds with the reversed inequality.
\ep

The proof of Proposition~\ref{P:compar} needs some preparations. If $i_1,\ldots,i_n$ are $n$ distinct elements of $S$ and $x\in\{0,1\}^S$, then we let $L^x_{i_1,\ldots,i_n}:\{0,1\}^n\to\{0,1\}$ denote the function
\be
L^x_{i_1,\ldots,i_n}(z_1,\ldots,z_n):=L(x')
\quad\mbox{with}\quad
x'(i):=\left\{\ba{ll}
\dis z_k\quad&\dis\mbox{if }i=i_k,\ 1\leq k\leq n,\\[5pt]
\dis x(i)\quad&\dis\mbox{if }i\not\in\{i_1,\ldots,i_n\}.
\ea\right.
\ee
We let $f_\vee$ and $f_\wedge$ denote the functions from $\{0,1\}^2$ to $\{0,1\}$ defined as
\be
f_\vee(z_1,z_2):=z_1\vee z_2\quand f_\wedge(z_1,z_2):=z_1\wedge z_2\qquad(z\in\{0,1\}^2).
\ee
For each $B\sub T$, we let $\psi^{-1}(B):=\{i\in S:\psi(i)\in B\}$ denote the inverse image of $B$ under $\psi$.

\bl[Contraction of two points]
In addition to the assumptions of Proposition~\ref{P:compar}, assume\label{L:contract} that there exist two elements $i_1,i_2\in S$ and an element $j_0\in T$ such that $\psi^{-1}(\{j_0\})=\{i_1,i_2\}$ and $\big|\psi^{-1}(\{j\})\big|=1$ for all $j\neq j_0$. Then
\be\label{contract}
\P\big[L(X)=1\big]-\P\big[L(Y\circ\psi)=1\big]
=p(1-p)\big\{\P\big[L^X_{i_1,i_2}=f_\vee\big]-\P\big[L^X_{i_1,i_2}=f_\wedge\big]\big\}.
\ee
\el

\bpro
Let $f_-,f_+,f_1$, and $f_2$ denote the functions from $\{0,1\}^2$ to $\{0,1\}$ defined as
\be
f_-(z_1,z_2):=0,\quad f_+(z_1,z_2):=1,\quad f_1(z_1,z_2):=z_1,\quand f_2(z_1,z_2):=z_2
\ee
$(z\in\{0,1\}^2)$. Then it is easy to check that $f_-,f_+,f_1,f_2,f_\vee,f_\wedge$ are all monotone Boolean functions from $\{0,1\}^2\to\{0,1\}$. Using the fact that the collection $X$ is i.i.d.\ with intensity $p$, we see that
\bc
\P\big[L(X)=1\big]&=&\dis\P\big[L^X_{i_1,i_2}=f_+\big]\\[5pt]
&&\dis+p\P\big[L^X_{i_1,i_2}=f_1\big]+p\P\big[L^X_{i_1,i_2}=f_2\big]\\[5pt]
&&\dis+\big(1-(1-p)^2\big)\P\big[L^X_{i_1,i_2}=f_\vee\big]+p^2\P\big[L^X_{i_1,i_2}=f_\wedge\big].
\ec
Similarly
\bc
\P\big[L(Y\circ\psi)=1\big]&=&\dis\P\big[L^X_{i_1,i_2}=f_+\big]\\[5pt]
&&\dis+p\P\big[L^X_{i_1,i_2}=f_1\big]+p\P\big[L^X_{i_1,i_2}=f_2\big]\\[5pt]
&&\dis+p\P\big[L^X_{i_1,i_2}=f_\vee\big]+p\P\big[L^X_{i_1,i_2}=f_\wedge\big].
\ec
Subtracting yields (\ref{contract}).
\epro

\bpro[of Proposition~\ref{P:compar}]
It suffices to prove the statement for minimal one-sets. The statement for minimal zero-sets then follows by the symmetry between zeros and ones. Let us say that $\psi$ is a \emph{pair contraction} if $\psi$ is surjective and  $|T|=|S|-1$. We first prove the statement under the additional assumption that $\psi$ is a pair contraction.

If $\psi$ is a pair contraction, then there exist two elements $i_1,i_2\in S$ and an element $j_0\in T$ such that $\psi^{-1}(\{j_0\})=\{i_1,i_2\}$ and $\big|\psi^{-1}(\{j\})\big|=1$ for all $j\neq j_0$, so Lemma~\ref{L:contract} is applicable and (\ref{compar}) follows from (\ref{contract}) provided we show that $L^x_{i_1,i_2}\neq f_\wedge$ for all $x\in\{0,1\}^S$. Assume that $L^x_{i_1,i_2}=f_\wedge$ for some $x\in\{0,1\}^S$. Since $L^x_{i_1,i_2}$ only depends on the values of $x(i)$ for $i\not\in\{i_1,i_2\}$, we can without loss of generality assume that $x(i_1)=x(i_2)=1$. Define $x_1,x_2\in\{0,1\}^S$ by
\be
x_1(i):=\left\{\ba{ll}
0\quad&\mbox{if }i=i_1,\\
x(i)\quad&\mbox{otherwise,}
\ea\right.
\quand
x_2(i):=\left\{\ba{ll}
0\quad&\mbox{if }i=i_2,\\
x(i)\quad&\mbox{otherwise.}
\ea\right.
\ee
Then $L^x_{i_1,i_2}=f_\wedge$ implies that $L(x)=1$ while $L(x_1)=L(x_2)=0$. This is possible only if $\{i_1,i_2\}\sub A$ for some $A\in\Ai(L)$, which contradicts (\ref{minint}). This completes the proof in the special case that $\psi$ is a pair contraction.

Before we prove the general statement, we make some elementary observations. If $S,T$ are finite sets and $\psi:S\to T$ is a function, then we can define $\Psi:\{0,1\}^T\to\{0,1\}^S$ by $\Psi(x):=x\circ\psi$. With this notation, (\ref{compar}) takes the form
\be\label{compar2}
\P\big[L\circ\Psi(Y)=1\big]\leq\P\big[L(X)=1\big].
\ee
We claim that
\be\label{AinA}
\forall A'\in\Ai(L\circ\Psi)\ \exists A\in\Ai(L)\mbox{ s.t.\ }A'\sub\psi(A).
\ee
Indeed, for each $y\in\{0,1\}^T$, one has
\be\ba{l}
\dis\exists A'\in\Ai(L\circ\Psi)\mbox{ s.t.\ }y(j)=1\ \forall j\in A'\\[5pt]
\dis\quad\desd\quad L\circ\Psi(y)=1\quad\desd\quad L(y\circ\psi)=1\\[5pt]
\dis\quad\desd\quad\exists A\in\Ai(L)\mbox{ s.t.\ }y\big(\psi(i)\big)=1\ \forall i\in A\\[5pt]
\dis\quad\desd\quad\exists A\in\Ai(L)\mbox{ s.t.\ }y(j)=1\ \forall j\in\psi(A),
\ec
and the implication $\volgt$ here can only hold for all $y\in\{0,1\}^T$ if (\ref{AinA}) is satisfied.

We now prove that (\ref{minint}) implies (\ref{compar2}). We can without loss of generality assume that $\psi$ is surjective. Then we can find an integer $n\geq 0$, finite sets $S_0,\ldots,S_n$ with $S_0=S$ and $S_n=T$, and pair contractions $\phi_k:S_{k-1}\to S_k$ $(1\leq k\leq n)$, such that $\psi=\phi_n\circ\cdots\phi_1$. For each $0\leq k\leq n$, let $L_k:\{0,1\}^{S_k}\to\{0,1\}$ and $\psi_k:S\to S_k$ be defined as
\be
L_k:=L\circ\Psi_k
\quad\mbox{with}\quad
\psi_k:=\phi_k\circ\cdots\circ\phi_1
\quand\Psi_k(x):=x\circ\psi_k\qquad\big(x\in\{0,1\}^{S_k}\big),
\ee
where $\psi_0$ is the identity map. For each $0\leq k\leq n$, let $X_k=(X_k(i))_{i\in S_k}$ be an i.i.d.\ collection of Boolean random variables with intensity $q$. We will prove (\ref{compar2}) by showing that
\be\label{compar3}
\P\big[L\circ\Psi(Y)=1\big]=\P\big[L_n(X_n)=1\big]\leq\cdots\leq\P\big[L_0(X_0)=1\big]=\P\big[L(X)=1\big].
\ee
Since $L_k=L_{k-1}\circ\Phi_k$ $(1\leq k\leq n)$ with $\Phi_k(x):=x\circ\phi_k$, it suffices to show that
\be
\P\big[L_{k-1}\circ\Phi_k(X_k)=1\big]\leq\P\big[L_{k-1}(X_{k-1})=1\big]\qquad(1\leq k\leq n).
\ee
Since $\phi_k$ is a pair contraction, by what we have already proved, it suffices to show that
\be\label{Apres}
\big|\phi_k(A')\big|=\big|A'\big|\qquad\forall A'\in\Ai(L_{k-1}).
\ee
The assumption (\ref{minint}) implies that $|\psi_n(A)|=|\psi(A)|=|A|=|\psi_0(A)|$ and hence $|\psi_0(A)|=\cdots=|\psi_n(A)|$ for all $A\in\Ai(L)$. In view of (\ref{AinA}), for each $A'\in\Ai(L_{k-1})$, there exists an $A\in\Ai(L)$ such that $A'\sub\psi_{k-1}(A)$. The fact that $|\phi_k\circ\psi_{k-1}(A)|=|\psi_k(A)|=|\psi_{k-1}(A)|$ now implies (\ref{Apres}) and the proof is complete.
\epro

\subsection{Comparison with the game on a tree}\label{S:treecomp}

In this subsection we prove Proposition~\ref{P:treecomp}. We will apply Proposition~\ref{P:compar} to the monotone Boolean function $L^{\rm AB}_{2n}$ from (\ref{FAB}). We start with some general observations.

Let $(G,\vec E,0,\tau)$ be a game-graph and let $L:\pa G\to\{0,1\}$ be the monotone Boolean map defined in (\ref{Fdef}). Using notation from Subsection~\ref{S:minmax}, let $\Si_\Ali$ and $\Si_\Bob$ denote the set of strategies of Alice and Bob, respectively, and as in (\ref{out}) let $o(\sig_\Ali,\sig_\Bob)$ denote the outcome of the game if Alice plays strategy $\sig_\Ali$ and Bob plays strategy $\sig_\Bob$. For each $\sig_\Ali\in\Si_\Ali$ and $\sig_\Bob\in\Si_\Bob$, we let
\be\label{AZstrat}
Z(\sig_\Ali):=\big\{o(\sig_\Ali,\sig'):\sig'\in\Si_\Bob\big\}
\quand
A(\sig_\Bob):=\big\{o(\sig',\sig_\Bob):\sig'\in\Si_\Ali\big\}
\ee
denote the sets of possible outcomes of the game if Alice plays the strategy $\sig_\Ali$ or Bob plays the strategy $\sig_\Bob$, respectively. We make the following simple observation.

\bl[Strategies and zero-sets]
Let\label{L:ZZZ}
\be\label{Ztrat}
\Zi_{\rm strat}:=\big\{Z(\sig_\Ali):\sig_\Ali\in\Si_\Ali\big\}
\quand
\Ai_{\rm strat}:=\big\{A(\sig_\Bob):\sig_\Bob\in\Si_\Bob\big\}.
\ee
Then
\be\label{ZZZ}
\Zi(L)\sub\Zi_{\rm strat}\sub\Zi^\up(L)
\quand
\Ai(L)\sub\Ai_{\rm strat}\sub\Ai^\up(L).
\ee
\el

\bpro
By symmetry, it suffices to prove the statement for zero-sets. We claim that
\be\label{LZ}
\Zi^\up(L)=\big\{Z\sub\pa G:\exists Z'\in\Zi_{\rm strat}\mbox{ s.t.\ }Z'\sub Z\big\}.
\ee
Indeed, using the fact that $L(x)=0$ if and only if Alice has a winning strategy given $x$, we see that
\[\ba{l}
\dis Z\in\Zi^\up(L)
\quad\desd\quad
L(x)=0\mbox{ with }x(v):=\left\{\ba{ll}
0\quad&\mbox{if }v\in Z\\
1\quad&\mbox{if }v\not\in Z
\ea\right.\\
\dis\quad\desd\quad
\mbox{Alice has a winning strategy given $x$}
\quad\desd\quad
\exists\sig_1\in\Si_1\mbox{ s.t.\ }x(v)=0\ \forall v\in Z(\sig_1)\\
\quad\desd\quad
\exists Z'\in\Zi_{\rm strat}\mbox{ s.t.\ }x(v)=0\ \forall x\in Z'
\quad\desd\quad
\exists Z'\in\Zi_{\rm strat}\mbox{ s.t.\ }Z'\sub Z,
\ea\]
which proves (\ref{LZ}). As an immediate consequence, we obtain
\[
Z\in\Zi_{\rm strat}
\quad\volgt\quad
Z':=Z\mbox{ satisfies }Z'\in\Zi_{\rm strat}\mbox{ and }Z'\sub Z
\quad\volgt\quad
Z\in\Zi^\up(L),
\]
i.e., $\Zi_{\rm strat}\sub\Zi^\up(L)$.

To prove that $\Zi(L)\sub\Zi_{\rm strat}$, assume that $Z\in\Zi(L)$. Then, since $\Zi(L)\sub\Zi^\up(L)$  (\ref{LZ}) implies that there exists a $Z'\in\Zi_{\rm strat}$ such that $Z'\sub Z$. Since $\Zi_{\rm strat}\sub\Zi^\up(L)$ we have $Z'\in\Zi^\up(L)$, so the minimality of $Z$ implies $Z'=Z$, proving that $Z\in\Zi_{\rm strat}$.
\epro

\bl[Projection property]
Let\label{L:project} $(D,\vec F,0)$ and $(D',\vec F',0')$ be decision graphs and assume that $\ell:D\to D'$ satisfies
\be\label{ell}
\ell(0)=0'\quand\mbox{$\ell$ maps $\Oi(v)$ surjectively into $\Oi(\ell(v))$ for all $v\in D$}.
\ee
Let $n\geq 1$ and let $L_n$ and $L'_n$ be the Boolean functions defined as in (\ref{Fdef}) in terms of the game-graphs $G_n(D)$ and $G_n(D')$, respectively. Then for each $n\geq 1$, $\ell$ maps $\pa G_n(D)$ into $\pa G_n(D')$, and
\be
L'_n(x)=L_n(x\circ\ell)\qquad\big(x\in\{0,1\}^{\pa G_n(D')}\big),
\ee
where $x\circ\ell\in\{0,1\}^{\pa G_n(D)}$ is defined as $x\circ\ell(v):=x\big(\ell(v)\big)$ $\big(v\in\pa G_n(D)\big)$.
\el

\bpro
It is easy to see that (\ref{ell}) implies that $\big|\ell(v)\big|=|v|$ $(v\in D)$. As a result, $\ell$ maps $\pa G_n(D)=\pa_nD$ into $\pa G_n(D')=\pa_nD'$. We can inductively extend the function $x:\pa_nD'\to\{0,1\}$ to a function $\ov x:[D']_n\to\{0,1\}$ such that for all $v\in\li D'\re_n$,
\be
\ov x(v)=\left\{\ba{ll}
\dis\bigwedge_{w\in\Oi(v)}\ov x(w)\quad&\mbox{if $|v|$ is even,}\\[5pt]
\dis\bigvee_{w\in\Oi(v)}\ov x(w)\quad&\mbox{if $|v|$ is odd.}
\ea\right.
\ee
Then $L'_n(x)=\ov x(0')$. Define $\ov y:[D]_n\to\{0,1\}$ by $\ov y(v):=\ov x\big(\ell(v)\big)$ $(v\in[D]_n)$. Then as a result of (\ref{ell}), we see that $\ov y$ satisfies the inductive relation
\be
\ov y(v)=\left\{\ba{ll}
\dis\bigwedge_{w\in\Oi(v)}\ov y(w)\quad&\mbox{if $|v|$ is even,}\\[5pt]
\dis\bigvee_{w\in\Oi(v)}\ov y(w)\quad&\mbox{if $|v|$ is odd.}
\ea\right.
\ee
Letting $y$ denote the restriction of $\ov y$ to $\pa_nD$, it follows that
\be
L_n(x\circ\ell)=L_n(y)=\ov y(0)=\ov x\big(\ell(0)\big)=\ov x(0')=L'_n(x).
\ee
\epro

An example of a function satisfying (\ref{ell}) is the function $\ell:\T\to\N^2$ defined as
\be\label{ellT}
\ell(\ibf):=\Big(\sum_{k=1}^n1_{\{i_k=1\}},\sum_{k=1}^n1_{\{i_k=2\}}\Big)\qquad(\ibf\in\T).
\ee
Further examples are the functions $\ell_1:\T\ltimes\T\to\N^2\ltimes\T$ and $\ell_2:\T\ltimes\T\to\T\ltimes\N^2$ defined as
\be\label{ell12}
\ell_1(\ibf,\jbf):=\big(\ell(\ibf),\jbf\big)
\quand
\ell_2(\ibf,\jbf):=\big(\ibf,\ell(\jbf)\big)
\qquad\big((\ibf,\jbf)\in\T\ltimes\T\big),
\ee
where $\ell$ is as in (\ref{ellT}). Lemma~\ref{L:project} now tells us that
\be\ba{r@{\,}c@{\,}ll}\label{Lproj}
\dis L^{\rm aB}_n(x)&=&\dis L^{\rm AB}_n(x\circ\ell_1)\quad&\dis\big(x\in\pa G_n(\N^2\ltimes\T)\big),\\[5pt]
\dis L^{\rm Ab}_n(x)&=&\dis L^{\rm AB}_n(x\circ\ell_2)\quad&\dis\big(x\in\pa G_n(\T\ltimes\N^2)\big).
\ec
We wish to apply Proposition~\ref{P:compar} to $L=L^{\rm AB}_{2n}$ and $\psi=\ell_1$ or $=\ell_2$. The following special property of the game-graph $G_n(\T\ltimes\T)$ will be crucial in our proof of Proposition~\ref{P:treecomp}.

\bl[The game-graph that is a tree]
For\label{L:funcstrat} the game-graph $G_n(\T\ltimes\T)$, let $\sig_\Ali\in\Si_\Ali$ and $\sig_\Bob\in\Si_\Bob$ be strategies for Alice and Bob, respectively, and let $Z(\sig_\Ali)$ and $A(\sig_\Bob)$ be defined as in (\ref{AZstrat}). Then for each $\jbf\in\pa_n\T$, there exists precisely one $\ibf\in\pa_n\T$ such that $(\ibf,\jbf)\in Z(\sig_\Ali)$. Likewise, for each $\ibf\in\pa_n\T$, there exists precisely one $\jbf\in\pa_n\T$ such that $(\ibf,\jbf)\in A(\sig_\Bob)$.
\el

\bpro
To prove the first statement, let $\sig_\Ali\in\Si_\Ali$ be a strategy for Alice and let $\jbf\in\pa_n\T$. Because of the tree structure of $\T$, given $\jbf$, we know exactly which moves Bob must make in each of his turns for the outcome of the game to be of the form $(\ibf,\jbf)$ for some $\ibf\in\pa_n\T$. Since we also know Alice's strategy, this means that we have all necessary information to determine the outcome of the game. This shows that $\ibf$ is unique. On the other hand, it is easy to see that there exists at least one strategy for Bob that leads to the outcome $(\ibf,\jbf)$. This completes the proof of the first statement. The second statement follows from the same argument.
\epro

\bpro[of Proposition~\ref{P:treecomp}]
Let
\be
X=\big(X(a,b)\big)_{a\in\pa_n\T,\ b\in\pa_n\T}
\quand
Y=\big(Y(a,b)\big)_{a\in\pa_n\T,\ b\in\pa_n\N^2}
\ee
be i.i.d.\ collections of Boolean random variables with intensity $p\in[0,1]$. We claim that
\be
P^{\rm Ab}_n(p)=\P\big[L^{\rm Ab}_{2n}(Y)=1\big]=\P\big[L^{\rm AB}_{2n}(Y\circ\ell_2)=1\big]\leq\P\big[L^{\rm AB}_{2n}(X)=1\big]=P^{\rm AB}_n(p).
\ee
Indeed, the first and final equalities follow from formula (\ref{Pminmax}), the second equality follows from formula (\ref{Lproj}), and the inequality will follow from Proposition~\ref{P:compar} provided we show that
\be\label{prescard}
\big|\ell_2(A)\big|=|A|\qquad\forall A\in\Ai(L^{\rm AB}_{2n}).
\ee
We claim that (\ref{prescard}) follows from Lemmas \ref{L:ZZZ} and \ref{L:funcstrat}. Indeed, Lemma~\ref{L:funcstrat} implies that each $A\in \Ai_{\rm strat}(L^{\rm AB}_{2n})$ has the property that for each $\ibf\in\pa_n\T$, there exists precisely one $\jbf\in\pa_n\T$ such that $(\ibf,\jbf)\in A$, so by Lemma~\ref{L:ZZZ} the same is true for  each $A\in\Ai(L^{\rm AB}_{2n})\sub\Ai_{\rm strat}(L^{\rm AB}_{2n})$. Since $\ell_2$ acts only on the second coordinate, this implies  (\ref{prescard}). This completes the proof that $P^{\rm Ab}_n(p)\leq P^{\rm AB}_n(p)$. The proof that $P^{\rm aB}_n(p)\geq P^{\rm AB}_n(p)$ follows from the same argument.
\epro

\subsection{Toom cycles}\label{S:cyc}

In Subsection~\ref{S:Pei}, we will derive upper bounds on the probability that Alice has a winning strategy in the games ${\rm Ab}_n(p)$ and ${\rm ab}_n(p)$, and on the probability that Bob has a winning strategy in the games ${\rm aB}_n(p)$ and ${\rm ab}_n(p)$, which then translate into some of the bounds in Propositions \ref{P:bounds} and \ref{P:ab}. To have a unified set-up, it will be convenient to introduce games ${\rm Ab}'_n(p)$ and ${\rm ab}'_n(p)$ that are identical to the games ${\rm Ab}_n(p)$ and ${\rm ab}_n(p)$ except that instead of Alice, Bob starts. Then the probability that Bob has a winning strategy in the games ${\rm aB}_n(p)$ and ${\rm ab}_n(p)$ is equal to the probability that Alice has a winning strategy in the games ${\rm Ab}'_n(1-p)$ and ${\rm ab}'_n(1-p)$, respectively. Therefore, it suffices to derive upper bounds on the probability that Alice has a winning strategy in the games ${\rm Ab}_n(p)$, ${\rm ab}_n(p)$, ${\rm Ab}'_n(p)$, and ${\rm ab}'_n(p)$.

We will use a Peierls argument first invented by Toom \cite{Too80} and further developed in \cite{SST22}. To prepare for this, in the present subsection, we formulate a theorem that says that if Alice has a winning strategy in any of these games, then a certain structure must be present in the game-graph that following \cite{SST22} we will call a ``Toom cycle''. We first need to introduce notation for the game-graphs of the games ${\rm Ab}'_n(p)$ and ${\rm ab}'_n(p)$.

Given a decision graph $D$ and integer $n\geq 1$, we define a game-graph
\be\label{Gacn}
G'_n(D)=(G'_n,\vec E'_n,0,\tau')
\ee
exactly as in (\ref{GnDD2a}), but with (\ref{GnDD2b}) replaced by
\be
\tau'(v):=\left\{\ba{ll}
\dis\Ali\quad&\dis\mbox{if }|v|\mbox{ is odd,}\\[5pt]
\dis\Bob\quad&\dis\mbox{if }|v|\mbox{ is even.}
\ea\right.\qquad(v\in\ring G_n).
\ee
In (\ref{ltimes}) and (\ref{ltimesF}), we defined a sort of ``product'' $D^1\ltimes D^2$ of two decision graphs $D^1$ and $D^2$. Elements of $D^1\ltimes D^2$ are pairs $(a,b)$ with $a\in D^1$ and $b\in D^2$. If $(a,b)$ is at distance $n$ from the root, then to reach a new state at distance $n+1$ from the root, we replace either $a$ or $b$ by a new state $a'$ or $b'$ in $D^1$ or $D^2$ that lies one step further from the root. Here we replace $a$ in case $n$ is even and $b$ in case $n$ is odd, as is natural for a game where Alice starts. For games where Bob starts, very much in the same spirit, we can define a second ``product'' $D^1\rtimes D^2$ which differs from $D^1\ltimes D^2$ only in the fact that we replace $a$ in case $n$ is odd and $b$ in case $n$ is even. Formally, we set
\be\label{rtimes}
D^1\rtimes D^2:=\big\{(a,b):(b,a)\in D^2\ltimes D^1\big\},
\ee
which is then equipped with the structure of a decision graph in the obvious way. We will be interested in the game-graphs
\be\label{Abgr}
G_n(\T\ltimes\N^2),\quad G'_n(\T\rtimes\N^2),\quad G_n(\N^2\ltimes\N^2),\quand G'_n(\N^2\rtimes\N^2).
\ee
Note that in each case, the decision graph of Bob is $\N^2$. We do not insist that $n$ is even, so we also allow for games where the same player has both the first and the last move. We will show that given a function that assigns winners to possible outcomes, if Alice has a winning strategy, then a ``Toom cycle'' must be present in these game-graphs.

\begin{figure}
\begin{center}
\inputtikz{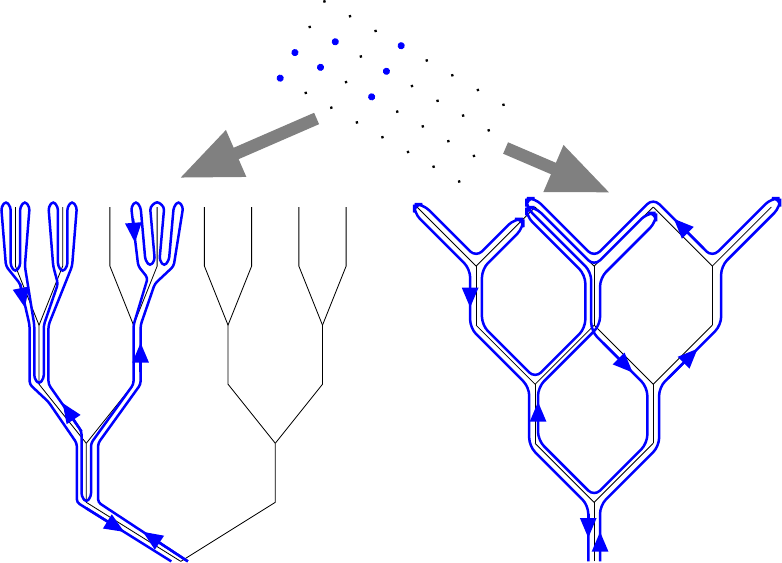}
\end{center}
\caption{A Toom cycle in the game-graph $G_6(\T\ltimes\N^2)$. The pictures on the left and right show two projections of the same Toom cycle (compare Figure~\ref{fig:gamegraph}). The picture above shows the possible outcomes $(a,b)$ of the game that the Toom cycle passes through. If Alice has a winning strategy, then there must exist a Toom cycle for which all these points correspond to a win for Alice.}\label{fig:Toom}
\end{figure}

Fix $n\geq 1$ and let $(G,\vec E,0,\tau)$ be any of the game-graphs in (\ref{Abgr}). Elements of $G$ are pairs $(a,b)$ with $a$ an element of $\T$ or $\N^2$ and $b$ an element of $\N^2$. Recall that if $\ibf=i_1\cdots i_m\in\T$ and $k\in\{1,2\}$, then $\ibf k\in\T$ denotes the word obtained by appending the letter $k$ to the word $\ibf$. To have a unified notation, it will be convenient to introduce for the decision-graph $\N^2$ the notation
\be\label{moves}
(i,j)k:=\left\{\ba{ll}
(i+1,j)\quad&\mbox{if }k=1,\\
(i,j+1)\quad&\mbox{if }k=2.
\ea\right.
\ee
Recall that $\ring{G}^\Ali$ and $\ring{G}^\Bob$ are the sets of internal states when it is Alice's and Bob's turn, respectively. Using notation as in (\ref{moves}), we set for $k\in\{1,2\}$
\be
\vec A_k:=\big\{\big((a,b),(ak,b)\big):(a,b)\in\ring{G}^\Ali\big\}
\quand
\vec B_k:=\big\{\big((a,b),(a,bk)\big):(a,b)\in\ring{G}^\Bob\big\}
\ee
and write $\vec A:=\vec A_1\cup\vec A_2$ and $\vec B:=\vec B_1\cup\vec B_2$ so that $\vec E=\vec A\cup\vec B$. For any set of directed edges $\vec F$, we let
\be
\lvec F:=\big\{(w,v):(v,w)\in\vec F\big\}
\ee
denote the set of directed edges obtained by reversing the direction of all edges in $\vec F$. A \emph{walk} in $G$ of \emph{length} $l$ is a finite word $\psi=\psi_0\cdots\psi_l$ made from the alphabet $G$ such that $(\psi_{k-1},\psi_k)\in\Vec E\cup\lvec E$ for each $0<k\leq l$. We call the substring $\psi_{k-1}\psi_k$ the \emph{$k$-th step} of the walk $\psi$ and we say that the $k$-th step is:
\begin{itemize}
\item[] \emph{straight-up} (su) if $(\psi_{k-1},\psi_k)\in\rvec A$, \emph{straight-down} (sd) if $(\psi_{k-1},\psi_k)\in\lvec A$,
\item[] \emph{left-up} (lu) if $(\psi_{k-1},\psi_k)\in\rvec B_2$, \emph{left-down} (ld) if $(\psi_{k-1},\psi_k)\in\lvec B_1$,
\item[] \emph{right-up} (ru) if $(\psi_{k-1},\psi_k)\in\rvec B_1$, \emph{right-down} (rd) if $(\psi_{k-1},\psi_k)\in\lvec B_2$.
\end{itemize}
This terminology is inspired by the right picture in Figure~\ref{fig:Toom}. Steps of the types su, ru, and lu are called \emph{up} and steps of the types sd, ld, and rd are called \emph{down}. Steps of the types su and sd are called \emph{straight}, those of the types lu and ld \emph{left}, and those of the types ru and rd \emph{right}. For any walk $\psi$, we set $N(\psi)=N:=\{0,\ldots,l\}$ and we partition $N$ into subsets $N_\up,N_\circ,N_\down,N_\ast$ defined as
\bc\label{NNNN}
\dis N_\up&:=&\dis\big\{k:0<k<l,\ (\psi_{k-1},\psi_k)\mbox{ is up},\ (\psi_k,\psi_{k+1})\mbox{ is up}\big\}\cup\{0\},\\[5pt]
\dis N_\circ&:=&\dis\big\{k:0<k<l,\ (\psi_{k-1},\psi_k)\mbox{ is down},\ (\psi_k,\psi_{k+1})\mbox{ is up}\big\},\\[5pt]
\dis N_\down&:=&\dis\big\{k:0<k<l,\ (\psi_{k-1},\psi_k)\mbox{ is down},\ (\psi_k,\psi_{k+1})\mbox{ is down}\big\}\cup\{l\},\\[5pt]
\dis N_\ast&:=&\dis\big\{k:0<k<l,\ (\psi_{k-1},\psi_k)\mbox{ is up},\ (\psi_k,\psi_{k+1})\mbox{ is down}\big\}.
\ec

We next define Toom cycles. A \emph{Toom cycle} in $G$ is a walk in $G$ of length $l\geq 2$ that must satisfy a number of requirements. First, we require that
\be
\psi_0=\psi_l=0,\mbox{ the first step is of type su or ru, and the last step is of type sd or rd.}
\ee
Next for $0<k<l$, depending on the type of the $k$-th step and on whether $\psi_k\in\ring{G}$ or $\in\pa G$, we will put constraints on the type of the $k+1$-th step. More precisely, if $\psi_k\in\ring{G}$, then only the following combinations are allowed for the types of the $k$-th and the $k+1$-th step:
\be\ba{lll}\label{steprules1}
\dis{\rm su}\to{\rm ru},\quad
&\dis{\rm ru}\to{\rm su},\quad
&\dis{\rm lu}\to{\rm su},\\[5pt]
\dis{\rm sd}\to{\rm rd}/{\rm ld},\quad
&\dis{\rm rd}\to{\rm sd},\quad
&\dis{\rm ld}\to{\rm lu}.
\ec
If, on the other hand, $\psi_k\in\pa G$, then only the following combinations are allowed:
\be\ba{lll}\label{steprules2}
\dis{\rm su}\to{\rm sd},\quad
&\dis{\rm ru}\to{\rm rd}/{\rm ld},\quad
&\dis{\rm lu}\to{\rm rd}/{\rm ld}.
\ec
See Figure~\ref{fig:Toom} for an illustration. Note that these rules imply that $N_\ast(\psi)=\{k:\psi_k\in\pa G\}$. In addition to the requirements above, a Toom cycle has to fulfill the following two additional, final requirements:
\begin{enumerate}
\item $\psi_k\neq\psi_m$ for all $k,m\in N_\ast$ with $k\neq m$,
\item if $\psi_k=\psi_m$ for some $k\in N_s$ and $m\in N_t$ with $s,t\in\{\up,\circ,\down\}$ and $s\leq t$ with respect to the total order on $\{\up,\circ,\down\}$ defined by $\up\,<\circ<\,\down$, then $k\leq m$.
\end{enumerate}
Requirement~(i) says that vertices in $\pa G$ are visited only once. Applying (ii) twice we see that if $\psi_k=\psi_m$ for some $k,m\in N_s$ with $s\in\{\up,\circ,\down\}$, then $k=m$, so (ii) can roughly be described by saying that internal vertices can be visited at most three times, and these visits have to take place in the right order depending on whether the walk is going up, is in a local minimum, or goes down.

Given a function $x:\pa G\to\{0,1\}$ that assigns a winner to each possible outcome, we say that a Toom cycle $\psi=\psi_0\cdots\psi_l$ is \emph{present} in $(G,x)$ if
\be
x(\psi_k)=0\quad\mbox{for all $0\leq k\leq l$ such that }\psi_k\in\pa G.
\ee
We will prove the following theorem.

\bt[Toom cycles]
Let\label{T:Toom} $G$ be any of the game-graphs in (\ref{Abgr}) and let $x:\pa G\to\{0,1\}$ be a function that assigns a winner to each possible outcome. If Alice has a winning strategy given $x$, then a Toom cycle is present in $(G,x)$.
\et

\begin{figure}
\begin{center}
\inputtikz{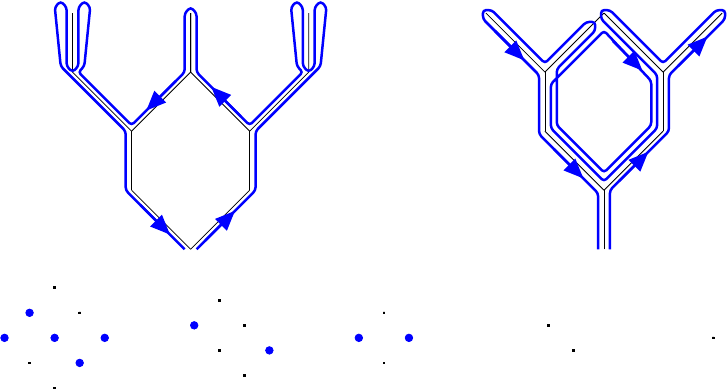}
\end{center}
\caption{A Toom cycle in the game-graph $G_4(\N^2\ltimes\N^2)$. The presence of this Toom cycle does not imply that Alice has a winning strategy. The pictures below show the zeros of the function $\ov x$ from (\ref{minmax}) at different levels of the graph, if at the top level all possible outcomes that the Toom cycle passes through have the value zero.}\label{fig:count}
\end{figure}

One can check that in Theorem~\ref{T:Toom}, the converse implication does not hold, i.e., the presence of a Toom cycle does not imply that Alice has a winning strategy. A counterexample for the game-graph $G_4(\N^2\ltimes\N^2)$ is shown in Figure~\ref{fig:count}. Intriguingly, we have not been able to construct a counterexample for the game-graphs $G_n(\T\ltimes\N^2)$ and $G'_n(\T\rtimes\N^2)$, leaving open the possibility that the converse implication holds in these cases. Theorem~\ref{T:Toom} is very similar to \cite[Thm~9]{SST22} and so is its proof (see in particular \cite[Figure~4]{SST22} for the idea of loop erasion), but since our setting is different and we will need some special properties of Toom cycles that are specific to our setting we will give an independent proof here. Strategies for Alice have an analogue for general monotone cellular automata. In this general context, they correspond to the so-called ``minimal explanations'' of \cite[Section~7.1]{SST22}, whose relation to Toom cycles is discussed in more detail in that article. We prove Theorem~\ref{T:Toom} in Subsection~\ref{S:Toomproof} and then apply it in Subsection~\ref{S:Pei}.

\subsection{Presence of Toom cycles}\label{S:Toomproof}

In this subsection we prove Theorem~\ref{T:Toom}. Let $G$ be any of the game-graphs in (\ref{Abgr}) and let $L$ be defined as in (\ref{Fdef}), i.e., this is the Boolean function defined by the requirement that if $x:\pa G\to\{0,1\}$ be a function that assigns a winner to each possible outcome, then $L(x)=0$ if and only if Alice has a winning strategy given $x$. Using notation introduced in Subsection~\ref{S:bolineq}, let $\Zi^\up(L)$ denote the set of zero-sets of $L$ and let $\Zi(L)$ denote its set of minimal zero-sets. Defining $\Zi_{\rm strat}$ as in (\ref{Ztrat}), we recall from Lemma~\ref{L:ZZZ} that
\be\label{ZZZ2}
\Zi(L)\sub\Zi_{\rm strat}\sub\Zi^\up(L).
\ee

\noi
\textbf{Remark} One can check that both inclusions in (\ref{ZZZ2}) are strict. We say that a strategy $\sig_\Ali$ for Alice is \emph{minimal} if the set $Z(\sig_\Ali)$ from (\ref{AZstrat}) satisfies $Z(\sig_\Ali)\in\Zi(L)$. Note that regardless of how one assigns winners to possible outcomes, it actually never makes sense for Alice to play a strategy that is not minimal, since for each non-minimal strategy she has at her disposal another strategy that is guaranteed to be at least as good in each situation. For the game-graphs $G_n(\T\ltimes\N^2)$ and $G'_n(\T\rtimes\N^2)$, it seems that with some effort, it is possible to give a precise description of all minimal strategies for Alice, but for the game-graphs $G_n(\N^2\ltimes\N^2)$ and $G'_n(\N^2\rtimes\N^2)$, classifying the minimal strategies seems a hard task.\med

\bpro[of Theorem~\ref{T:Toom}]
For each Toom cycle $\psi=\psi_0\cdots\psi_l$, we write
\be
Z_\psi:=\big\{\psi_k:0<k<l,\ \psi_k\in\pa G\big\},
\ee
and we set
\be\label{ZToom}
\Zi_{\rm Toom}:=\big\{Z_\psi:\psi\mbox{ is a Toom cycle}\big\}.
\ee
We will prove that
\be\label{nec}
\forall Z\in\Zi_{\rm strat}\ \exists Z'\in\Zi_{\rm Toom}\mbox{ s.t.\ }Z'\sub Z,
\ee
which clearly implies the statement of the theorem. The proof is by induction on $n$. For concreteness, we write down the proof for the game-graphs $G_n(\T\ltimes\N^2)$ and $G_n(\N^2\ltimes\N^2)$, and at the end remark on how the proof has to be adapted for the other two game-graphs. To have a unified notation, we write $D\ltimes\N^2$ where $D=\T$ or $=\N^2$. Let $\Zi_{\rm strat}(n)$ and $\Zi_{\rm Toom}(n)$ denote the sets defined in (\ref{Ztrat}) and (\ref{ZToom}) for a given value of $n$. It will be convenient to also include the case $n=0$. For this purpose, we use the convention that in this case there is precise one Toom cycle, which is the walk of length $l=0$ given by $\psi=\psi_0=(0,0)$, where $(0,0)$ denotes the root of $D\ltimes\N^2$. Then $\Zi_{\rm strat}(0)$ and $\Zi_{\rm Toom}(0)$ both have a single element, which is the singleton $\{(0,0)\}$, and (\ref{nec}) holds for $n=0$.

If $n$ is even, $Z\in\Zi_{\rm strat}(n)$ is given, and $Z\ni(a,b)\mapsto\kappa_{(a,b)}\in\{1,2\}$ is a function, then we can define $\ti Z\in\Zi_{\rm strat}(n+1)$ by
\be\label{tiZeven}
\ti Z:=\big\{(a\kappa_{(a,b)},b):(a,b)\in Z\big\},
\ee
where $a\kappa_{(a,b)}$ is defined as in (\ref{moves}) if the decision graph of Alice is $\N^2$. Formula (\ref{tiZeven}) corresponds to Alice playing the strategy that defined $Z$ up to the $n$-th turn of the game and then playing the move $\kappa_{(a,b)}$ if at that moment the state of the game is $(a,b)$. We see from this that $\ti Z\in\Zi_{\rm strat}(n+1)$, and each element of $\Zi_{\rm strat}(n+1)$ is of this form for some $Z\in\Zi_{\rm strat}(n)$. Similarly, if $n$ is odd and $Z\in\Zi_{\rm strat}(n)$ is given, then we can define $\ti Z\in\Zi_{\rm strat}(n+1)$ by
\be\label{tiZodd}
\ti Z:=\big\{(a,b1),(a,b2):(a,b)\in Z\big\}.
\ee
Indeed, this corresponds to Alice playing the strategy that defined $Z$ up to the $n$-th turn of the game, while Bob has now one extra final move. We see from this that $\ti Z\in\Zi_{\rm strat}(n+1)$, and each element of $\Zi_{\rm strat}(n+1)$ is of this form for some $Z\in\Zi_{\rm strat}(n)$. In view of this, to complete the induction step of the proof, it suffices to show that:
\begin{itemize}
\item[] If for $Z\in\Zi_{\rm strat}(n)$ there exists a $Z'\in\Zi_{\rm Toom}(n)$ such that $Z'\sub Z$ and $\ti Z$ is defined as in (\ref{tiZeven}) or (\ref{tiZodd}) depending on whether $n$ is even or odd, then there exists a $\ti Z'\in\Zi_{\rm Toom}(n+1)$ such that $\ti Z'\sub\ti Z$.
\end{itemize}

To prove this, let $\psi=\psi_0\cdots\psi_l$ be a Toom cycle such that $Z'=Z_\psi$. If $n$ is even, then we modify $\psi$ in such a way that for each $k$ such that $\psi_k\in\pa G_n(D\ltimes\N^2)$, in place of $\psi_k=(a,b)$ we insert the string $(a,b)(a\kappa_{(a,b)},b)(a,b)$. Note that this adds two steps to the Toom cycle, the first of which is straight-up and the second straight-down. If $n$ is odd, then we modify $\psi$ in such a way that for each $k$ such that $\psi_k\in\pa G_n(D\ltimes\N^2)$, in place of $\psi_k=(a,b)$ we insert the string $(a,b)(a,b1)(a,b)(a,b2)(a,b)$. Note that this adds four steps to the Toom cycle: right-up, left-down, left-up, and right-down.

Let $\ti\psi$ denote the modified walk. Then it is straightforward to check that $\ti Z':=Z_{\ti\psi}\sub\ti Z$ and $\ti\psi$ satisfies all requirements of a Toom cycle, except possibly for requirement~(i), i.e., it may now be the case that some potential outcomes of the game are visited more than once. To remedy this, we apply the procedure of \emph{loop erasion} which is also the main idea behind the proof of \cite[Thm~9]{SST22}. If $\ti\psi_k=\ti\psi_{k'}\in\pa G_{n+1}$ for some $k<k'$, then we remove the string $\ti\psi_{k+1}\cdots\ti\psi_{k'}$ from $\ti\psi$. One can check that this preserves all the requirements of a Toom cycle that $\ti\psi$ satisfies. Moreover, by applying this procedure till it is no longer possible to do so, we can make sure that also requirement~(i) becomes valid again. In this way, we find a new modified Toom cycle $\ti\psi$ such that $\ti Z':=Z_{\ti\psi}\sub\ti Z$. Note, however, that because of the loop erasion, even if we started with $Z'=Z$, it may happen that $\ti Z'$ is a strict subset of $\ti Z$. 

The proof for the game-graphs $G'_n(\T\rtimes\N^2)$ and $G'_n(\N^2\rtimes\N^2)$ is essentially the same. It differs only in the sense that the roles of even and odd $n$ are interchanged.
\epro

\noi
\textbf{Remark} For the game-graphs $G_n(\N^2\ltimes\N^2)$ and $G'_n(\N^2\rtimes\N^2)$, one can check that because of the loop erasion, the construction in the proof of Theorem~\ref{T:Toom} can yield Toom cycles in which a right-up step is followed by a right-down step, as explicitly allowed in (\ref{steprules2}). For the game-graphs $G_n(\T\ltimes\N^2)$ and $G'_n(\T\rtimes\N^2)$, it seems that on the other hand such transitions never arise from the construction in the proof of Theorem~\ref{T:Toom}, so for these game-graphs it should be possible to prove a stronger version of the theorem in which the Toom cycle satisfies additional properties. Since we do not need this for our main results, we do not pursue this further.

\subsection{The Peierls argument}\label{S:Pei}

Recall the definition of the games ${\rm Ab}'_n(p)$ and ${\rm ab}'_n(p)$ at the beginning of Subsection~\ref{S:cyc}. In this subsection, we apply Theorem~\ref{T:Toom} to derive upper bounds on the probability that Alice has a winning strategy in the games ${\rm Ab}_n(p)$, ${\rm ab}_n(p)$, ${\rm Ab}'_n(p)$, and ${\rm ab}'_n(p)$. We start by considering game-graphs where Alice starts, since in this case the argument is slightly easier.

\bl[Toom cycles if Alice starts]
If\label{L:step} $\psi$ is a Toom cycle in $G=G_{2n}(\T\ltimes\N^2)$ or $G_{2n}(\N^2\ltimes\N^2)$, then there exists an integer $m\geq 1$ such that the number of steps of each of the types su, sd, lu, ld, ru, rd equals $m$, and the Toom cycle passes through precisely $m+1$ different possible outcomes of the game.
\el

\bpro
We see from (\ref{steprules1}) and  (\ref{steprules2}) that after a step of a type in $\{{\rm su},{\rm ld},{\rm rd}\}$ there always follows a step of a type in $\{{\rm ru},{\rm lu},{\rm sd}\}$ and vice versa. Therefore, since the first step of the Toom cycle is straight-up and the last step is straight-down, the number of steps of the Toom cycle is even and if we divide its steps in pairs of two consecutive steps, then the first step of such a pair is always of a type in $\{{\rm su},{\rm ld},{\rm rd}\}$ and the second step of a type in $\{{\rm ru},{\rm lu},{\rm sd}\}$. It follows from the structure of $G=G_{2n}(\T\ltimes\N^2)$ and $G_{2n}(\N^2\ltimes\N^2)$ that a straight-up step cannot end in $\pa G_{2n}$, so by (\ref{steprules1}) only three types of pairs of two consecutive steps are possible:
\be\label{steppair}
{\rm su}\mbox{ followed by }{\rm ru},\quad
{\rm ld}\mbox{ followed by }{\rm lu},\quad\mbox{or}\quad
{\rm rd}\mbox{ followed by }{\rm sd}.
\ee
By a slight abuse of notation, let us also write ${\rm su}$ to denote the number of straight-up steps in the Toom cycle, and similarly for the other five types of steps. Then (\ref{steppair}) implies that
\be\label{paireq}
{\rm su}={\rm ru},\quad{\rm ld}={\rm lu},\quad{\rm rd}={\rm sd}.
\ee
Since the distance from the root increases by one in each up move and decreases by one in each down move, and the cycle ends where it began, the number of up moves must equal the number of down moves, which gives
\be\label{updo}
{\rm su}+{\rm lu}+{\rm ru}={\rm sd}+{\rm ld}+{\rm rd}.
\ee
Similarly, if $f\big(a,(i,j)\big):=i-j$ denotes the difference between number of 1 moves and 2 moves played by Bob when the state of the game is $(a,(i,j))$, then the function $f$ increases by one in each right move, decreases by one in each left move, and stays the same in straight moves, so we see that the number of right moves must equal the number of left moves, which gives
\be\label{lr}
{\rm lu}+{\rm ld}={\rm ru}+{\rm rd}.
\ee
Combining (\ref{paireq}), (\ref{updo}), and (\ref{lr}), we see that there exists an integer $m\geq 1$ such that
\be
{\rm su}={\rm sd}={\rm lu}={\rm ld}={\rm ru}={\rm rd}=:m.
\ee
In view of (\ref{steprules1}), each change from the down direction to the up direction is associated with a pair consisting of a left-down step followed by a left-up step, so the number of these changes is precisely $m$. Since the first step is up and the last step is down, we change from the up direction to the down direction once more than the other way around. Because of (\ref{steprules1}) and  (\ref{steprules2}), changes from the up direction to the down direction can only happen in states that are a possible outcome of the game, and because of the crucial condition~(i) in the definition of a Toom cycle, the walk visits each possible outcome at most once. It follows that the Toom cycle passes through precisely $m+1$ different possible outcomes of the game.
\epro

\bp[Peierls estimates if Alice starts]
For\label{P:APei} each $n\geq 1$, the probability that Alice has a winning strategy in the games ${\rm Ab}_n(p)$ and ${\rm ab}_n(p)$ can be estimated from above by
\be
1-P^{\rm Ab}_n(p)\leq\sum_{m=n}^\infty 2^{3m}(1-p)^{m+1}
\quand
1-P^{\rm ab}_n(p)\leq\sum_{m=n}^\infty 2^{4m}(1-p)^{m+1}.
\ee
In particular,
\be
P^{\rm Ab}_n(p)\asto{n}1\quad\mbox{for all }p>\frac{7}{8}
\qquad\mbox{and}\qquad
P^{\rm ab}_n(p)\asto{n}1\quad\mbox{for all }p>\frac{15}{16}.
\ee
\ep

\bpro
Let $G=G_{2n}(\T\ltimes\N^2)$ and let $\big(X(a,b)\big)_{(a,b)\in\pa G}$ be i.i.d.\ Boolean random variables with intensity $p$. By Theorem~\ref{T:Toom}, the probability that Alice has a winning strategy in the game ${\rm Ab}_n(p)$ can be estimated from above by the probability that a Toom cycle is present in $(G,X)$, which in turn can be estimated from above by the expected number of such Toom cycles. By Lemma~\ref{L:step}, the length of such a Toom cycle must be of the form $6m$ with $m$ an integer. We have $m\geq n$ since the cycle must make at least $n$ straight-up steps before it can start walking down. Since by Lemma~\ref{L:step}, a Toom cycle of length $6m$ visits precisely $m+1$ possible outcomes, the probability that a given Toom cycle of length $6m$ is present in $(G,X)$ is precisely $(1-p)^{m+1}$. This gives the bound
\be
1-P^{\rm Ab}_n(p)\leq\sum_{m=n}^\infty M_m(1-p)^{m+1},
\ee
where $M_m$ denotes the number of distinct Toom cycles of length $6m$ in $G$. We claim that
\be
M_m\leq 2^m\cdot 2^{m+1}\cdot 2^{m-1}.
\ee
Here the first factor comes from the fact that for each straight-up step we have two choices (corresponding to the two possible moves of Alice), the second factor comes from the fact that each time we arrive in a possible outcome of the game, we have to choose whether the next step is left-down or right-down, and the third factor comes from the fact that after each straight-down step except the very last one, we have to choose whether the next step is left-down or right-down. By (\ref{steprules1}) and (\ref{steprules2}) these are the only choices we need to make to uniquely determine a Toom cycle.

The argument for the game ${\rm ab}_n(p)$ is precisely the same, except that now the decision graph of Alice is no longer a tree which has the effect that also at the beginning of each straight-down step we may have to choose between two possible ways of making a straight-down step (see Figure~\ref{fig:count}).
\epro

For the game-graphs where Bob starts, the argument is only slightly more complicated.

\bp[Peierls estimates if Bob starts]
For\label{P:BPei} each $n\geq 1$, the probability that Alice has a winning strategy in the games ${\rm Ab}'_n(p)$ and ${\rm ab}'_n(p)$ can be estimated from above by
\be\label{BPei}
\sum_{m=n}^\infty 2^{4m}(1-p)^{m+1}
\quand
\sum_{m=n}^\infty 2^{6m}(1-p)^{m+1}.
\ee
As a result,
\be\label{Alwin}
P^{\rm aB}_n(p)\asto{n}0\quad\mbox{for all }p<\frac{1}{16}
\qquad\mbox{and}\qquad
P^{\rm ab}_n(p)\asto{n}0\quad\mbox{for all }p<\frac{1}{64}.
\ee
\ep

\bpro
The proof is very similar to the proof of Proposition~\ref{P:APei}, but we have to modify Lemma~\ref{L:step}. Toom cycles in the game-graphs $G'_{2n}(\T\rtimes\N^2)$ and $G'_{2n}(\N^2\rtimes\N^2)$ can be obtained from Toom cycles in the game-graphs $G_{2n}(\T\ltimes\N^2)$ and $G_{2n}(\N^2\ltimes\N^2)$ by removing the first and last steps (which were straight-up and straight-down), and then inserting at each instance where the Toom cycle changes from the up to the down direction one extra straight-up step followed my a straight-down step. The argument then proceeds as before up to the point where we count the number of modified Toom cycles for which the original (unmodified) Toom cycle had length $6m$. For the game-graph $G'_{2n}(\T\rtimes\N^2)$, we get a factor $2$ for each straight-up step. Due to the modification, the number of straight-up steps has decreased by one and increased by $m+1$, so we gain an extra factor $2^m$. For the game-graph $G'_{2n}(\N^2\rtimes\N^2)$, we also get a factor $2$ for each straight-down step, so we gain a factor $2^{2m}$. This explains (\ref{BPei}).

Formula (\ref{Alwin}) follows from (\ref{BPei}) and the fact that the probability that Bob has a winning strategy in the games ${\rm aB}_n(p)$ and ${\rm ab}_n(p)$ is equal to the probability that Alice has a winning strategy in the games ${\rm Ab}'_n(1-p)$ and ${\rm ab}'_n(1-p)$, respectively.
\epro

The bounds for $P^{\rm ab}_n(p)$ from Propositions \ref{P:APei} and \ref{P:BPei} immediately imply Proposition~\ref{P:ab}, while the bounds for $P^{\rm Ab}_n(p)$ and $P^{\rm aB}_n(p)$ will translate into the upper bound on $p^{\rm Ab}_{\rm c}$ and lower bound on $p^{\rm aB}_{\rm c}$ from Proposition~\ref{P:bounds}, once sharpness of the transition is established.

\subsection{Sharpness of the transition}

Throughout this subsection, $G_{2n}$ denotes any of the game-graphs $G_{2n}(\T\ltimes\N^2)$ and $G_{2n}(\N^2\ltimes\T)$, which are the game-graphs of the games ${\rm Ab}_n(p)$ and ${\rm aB}_n(p)$, respectively, and $L_{2n}$ denotes any of the Boolean functions $L^{\rm Ab}_{2n}$ and $L^{\rm aB}_{2n}$ that tell us who has a winning strategy for each given assignment $x=\big(x(v)\big)_{v\in\pa G_{2n}}$ of winners to possible outcomes, with $0$ and $1$ indicating a win for Alice and Bob, respectively.

For each $p\in[0,1]$, we let $\big(X^p(v)\big)_{v\in\pa G_{2n}}$ be i.i.d.\ Boolean variables with intensity $p$ and for each $v\in\pa G_{2n}$, we define
\be
X^p_{v,0}(w):=\left\{\ba{ll}
0\quad&\mbox{if }w=v,\\[5pt]
X^p(w)\quad&\mbox{otherwise,}
\ea\right.
\quand
X^p_{v,1}(w):=\left\{\ba{ll}
1\quad&\mbox{if }w=v,\\[5pt]
X^p(w)\quad&\mbox{otherwise,}
\ea\right.
\ee
$(w\in\pa G_{2n})$. By definition, a possible outcome $v\in\pa G_{2n}$ is \emph{pivotal} for $L_{2n}$ if $L_{2n}(X^p_{v,0})\neq L_{2n}(X^p_{v,1})$, which by the monotonicity of $L_{2n}$ means that $L_{2n}(X^p_{v,0})=0$ and $L_{2n}(X^p_{v,1})=1$. The \emph{influence} of $v$ is defined as
\be
I^p_{2n}(v):=\P\big[v\mbox{ is pivotal for }L_{2n}\big]\qquad(v\in\pa G_{2n})
\ee
and Russo's formula (see, e.g., \cite[Thm 3.2]{GS14}) tells us that the probability $P_n(p):=\P[L_{2n}(X^p)=1]$ that Bob has a winning strategy is continuously differentiable with
\be\label{eq:russo}
\dif{p}P_n(p)=\sum_{v\in\pa G_{2n}}I^p_{2n}(v).
\ee
It is easy to see that each possible outcome has a positive influence, so Russo's formula implies that $P_n$ is strictly increasing. Since moreover $P_n(0)=0$ and $P_n(1)=1$, we conclude that the function $P_n$ is invertible. We will prove the following result.

\bt[Sharpness of the transition]
For\label{T:sharp} each $0<\eps\leq\ha$, there exists a constant $C<\infty$ such that
\be\label{sharp}
P_n^{-1}(1-\eps)-P_n^{-1}(\eps)\leq\frac{C}{n}\qquad(n\geq 1).
\ee
\et

We start with a preparatory lemma.

\bl[Lower bound on the total influence]
There\label{L:inflow} exists a constant $c>0$ such that for all $p\in[0,1]$, one has
\be\label{inflow}
\sum_{v\in\pa G_{2n}}I^p_{2n}(v)\geq\inf_{J>0}\big[2^nJ\vee c\var\big(L_{2n}(X^p)\big)\log(1/J)\big].
\ee
\el

\bpro
Let us write
\be
I^p_{2n}:=\sum_{v\in\pa G_{2n}}I^p_{2n}(v)\quand J^p_{2n}:=\max_{v\in\pa G_{2n}}I^p_{2n}(v).
\ee
For the game-graph $G_{2n}=G_{2n}(\T\ltimes\N^2)$, the set of possible outcomes is $\pa G_{2n}=\pa_n\T\times\pa_n\N^2$, where $\pa_n\T$ has $2^n$ elements that play a completely symmetric role in $G_{2n}$, so that the influence $I^p_{2n}(a,b)$ of an element $(a,b)\in\pa G_{2n}$ depends on $b$ but not $a$. Choosing $b_\ast\in\pa_n\N^2$ such that $I^p_{2n}(a,b_\ast)=J^p_{2n}$ for all $a\in\pa_n\T$ gives the trivial bound
\be
I^p_{2n}\geq\sum_{a\in\pa_n\T}I^p_{2n}(a,b_\ast)=2^nJ^p_{2n}.
\ee
For the game-graph $G_{2n}=G_{2n}(\N^2\ltimes\T)$, we get the same bound by reversing the roles of $a$ and $b$. To complete the proof, we apply an inequality from the theory of noise sensitivity. By \cite[Thm~3.4]{GS14}, which is based on \cite{BK+92}, there exists a universal constant $c>0$ such that
\be
I^p_{2n}\geq c\var\big(L_{2n}(X^p)\big)\log(1/J^p_{2n}).
\ee
Combining these two estimates one arrives at (\ref{inflow}).
\epro

\bpro[of Theorem~\ref{T:sharp}]
Since $P_n:[0,1]\to[0,1]$ is strictly increasing we have
\be
\var\big(L_{2n}(X^p)\big)=P_n(p)\big(1-P_n(p)\big)
\geq\eps(1-\eps)\qquad\big(p\in[P_n^{-1}(\eps),P_n^{-1}(1-\eps)]\big).
\ee
Russo's formula (\ref{eq:russo}) and Lemma~\ref{L:inflow} now tell us that there exists a constant $c>0$ such that
\be\label{difineq}
\dif{p}P_n(p)\geq\inf_{J>0}\big[2^nJ\vee c\eps(1-\eps)\log(1/J)\big]\qquad\big(p\in[P_n^{-1}(\eps),P_n^{-1}(1-\eps)]\big).
\ee
Making $c$ smaller if necessary, we can without loss of generality assume that
\be\label{wlog}
c\eps(1-\eps)\log 2\leq 1.
\ee
To estimate the right-hand side of (\ref{difineq}) from below, we set $J_n:=n2^{-n}$. Since the function $J\mapsto 2^nJ$ is increasing and $J\mapsto c\eps(1-\eps)\log(1/J)$ is decreasing, we have
\be\label{lrest}
2^nJ\vee c\eps(1-\eps)\log(1/J)\left\{\ba{ll}
\geq 2^nJ_n\quad&\mbox{for }J\geq J_n,\\[5pt]
\geq c\eps(1-\eps)\log(1/J_n)\quad&\mbox{for }J\leq J_n.
\ea\right.
\ee
We observe that by (\ref{wlog}),
\be\label{event}
c\eps(1-\eps)\log(1/J_n)=c\eps(1-\eps)(n\log 2-\log n)\leq n=2^nJ_n.
\ee
Choose $0<c'<c\eps(1-\eps)\log 2$. Then (\ref{difineq}), (\ref{lrest}), and (\ref{wlog}) combine to give
\be
\dif{p} P_n(p)\geq c\eps(1-\eps)(n\log 2-\log n)\geq c'n\qquad\big(p\in[P_n^{-1}(\eps),P_n^{-1}(1-\eps)]\big).
\ee
Setting $C:=1/c'$, we obtain (\ref{sharp}).
\epro

\subsection{The cellular automata}

In this subsection we study the cellular automata $(X^{{\rm Ab},p}_t)_{t\geq 0}$ and $(X^{{\rm aB},p}_t)_{t\geq 0}$ from Subsection~\ref{S:cell}. It will be convenient to change the notation a bit. Given a random variable $X_0$ with values in $\{0,1\}^{\N^2}$, we can define a discrete-time process $(X_t)_{t\geq 0}$ with values in $\{0,1\}^{\N^2}$ by setting
\be\label{Xdef}
X_t(i,j):=\left\{\ba{ll}
\dis X_{t-1}(2i,j)\wedge X_{t-1}(2i+1,j)\quad&\mbox{if $t$ is even},\\[5pt]
\dis X_{t-1}(i,j)\vee X_{t-1}(i,j+1)\quad&\mbox{if $t$ is odd},
\ea\right.\qquad\big(t>0,\ (i,j)\in\N^2\big).
\ee
In particular, if $\big(X_0(i,j)\big)_{(i,j)\in\N^2}$ are i.i.d.\ with intensity $p$, then this is the cellular automaton $(X^{{\rm Ab},p}_t)_{t\geq 0}$ defined in Subsection~\ref{S:cell}. Similarly, given a random variable $\hX_0$ with values in $\{0,1\}^{\N^2}$, we can define a discrete-time process $(\hX_t)_{t\geq 0}$ with values in $\{0,1\}^{\N^2}$ by setting
\be\label{hatXdef}
\hX_t(i,j):=\left\{\ba{ll}
\dis\hX_{t-1}(2i,j)\wedge\hX_{t-1}(2i+1,j)\quad&\mbox{if $t$ is odd},\\[5pt]
\dis\hX_{t-1}(i,j)\vee\hX_{t-1}(i,j+1)\quad&\mbox{if $t$ is even},
\ea\right.\qquad\big(t>0,\ (i,j)\in\N^2\big),
\ee
where compared to (\ref{Xdef}) all we have done is that we have interchanged the rules at even and odd times. The cellular automaton $(\hX_t)_{t\geq 0}$ is similar to the cellular automaton $(X^{{\rm aB},p}_t)_{t\geq 0}$ from Subsection~\ref{S:cell}, but compared to the definitions there, the roles of $i$ and $j$ and the roles of $0$ and $1$ are interchanged. More precisely, if $\big(\hX_0(i,j)\big)_{(i,j)\in\N^2}$ are i.i.d.\ with intensity $1-p$, then setting
\be\label{hataB}
X^{{\rm aB},p}_t(i,j):=1-\hX_t(j,i)\qquad\big(t\geq 0,\ (i,j)\in\N^2\big)
\ee
yields a cellular automaton $(X^{{\rm aB},p}_t)_{t\geq 0}$ as defined in Subsection~\ref{S:cell}.

As observed in (\ref{column}), the cellular automaton $(X^{{\rm Ab},p}_t)_{t\geq 0}$ has the property that at each time its columns are independent. This holds not only for product initial laws but more generally if the initial state $X_0$ has the property that its columns are independent, and the same holds for $(\hX_t)_{t\geq 0}$. This motivates the following definitions.

We let $S:=\{0,1\}^\N$ denote the space of all functions $y:\N\to\{0,1\}$. We equip $S$ with the product topology and the associated Borel-\si-field and we set
\be\label{MiSi}
\Mi:=\big\{\mu:\mu\mbox{ is a probability law on }S\big\}
\quand
\Si:=\big\{\mu\in\Mi:\mu\mbox{ is stationary}\big\}.
\ee
We equip $\Mi$ with the topology of weak convergence of probability laws. Then $\Si$ is a closed subset of $\Mi$. We define maps $\Phi_{\rm A}:S^2\to S$ and $\Phi_{\rm b}:S\to S$ by
\be
\Phi_{\rm A}(y,z)(k):=y(k)\wedge z(k)
\quand
\Phi_{\rm b}(y)(k):=y(k)\vee y(k+1)\qquad(k\in\N),
\ee
we define maps $\Fi_{\rm A}$ and $\Fi_{\rm b}$ from $\Mi$ to $\Mi$ by
\be\ba{l}
\dis\Fi_{\rm A}(\mu):=\P\big[\Phi_{\rm A}(Y,Z)\in\,\cdot\,\big]
\quand
\Fi_{\rm b}(\mu):=\P\big[\Phi_{\rm b}(Y)\in\,\cdot\,\big]\\[5pt]
\qquad\mbox{where $Y$ and $Z$ are independent with law $\mu\in\Mi$,}
\ec
and we set
\be\label{Fi}
\Fi:=\Fi_{\rm A}\circ\Fi_{\rm b}\quand\hFi:=\Fi_{\rm b}\circ\Fi_{\rm A}.
\ee
The following lemma is a more precise formulation of the observation (\ref{column}).

\bl[Law of the columns]
Let\label{L:column} $(X_t)_{t\geq 0}$ be the cellular automaton defined in (\ref{Xdef}) started in an initial law such that
\be
\big(X_0(i,\,\cdot\,)\big)_{i\in\N}\mbox{ are i.i.d.\ with common law }\mu.
\ee
Then for each $n\in\N$,
\be\label{column2}
\big(X_{2n}(i,\,\cdot\,)\big)_{i\in\N}\mbox{ are i.i.d.\ with common law }\Fi^n(\mu),
\ee
where $\Fi^n$ denotes the $n$-th iterate of the map defined in (\ref{Fi}). The same is true with $(X_t)_{t\geq 0}$ replaced by $(\hX_t)_{t\geq 0}$ and $\Fi$ replaced by $\hFi$.
\el

\bpro
The proof is by induction on $n$. Assuming that the statement holds for $n-1$, we see from (\ref{Xdef}) that
\be
X_{2n-1}(i,\,\cdot\,)=\Phi_{\rm b}\big(X_{2n-2}(i,\,\cdot\,)\big)\qquad(i\in\N),
\ee
which implies (compare Figure~\ref{fig:cellaut}) that
\be
\big(X_{2n-1}(i,\,\cdot\,)\big)_{i\in\N}\mbox{ are i.i.d.\ with common law }\Fi_{\rm b}\circ\Fi^{n-1}(\mu).
\ee
Using again (\ref{Xdef}), we see that
\be
X_{2n}(i,\,\cdot\,)=\Phi_{\rm A}\big(X_{2n-1}(2i,\,\cdot\,),X_{2n-1}(2i+1,\,\cdot\,)\big)\qquad(i\in\N),
\ee
which implies (\ref{column2}). The proof for $(\hX_t)_{t\geq 0}$ is the same, but we apply $\Phi_{\rm A}$ and $\Phi_{\rm b}$ in the opposite order.
\epro

We let $\pi_p$ denote the product measure on $S=\{0,1\}^\N$ with intensity $p$ and for any $\mu\in\Si$, we let $\li\mu\re:=\int\mu(\di y)y(i)$ which by stationarity does not depend on $i\in\N$. Then our main functions of interest can be expressed in terms of the ``column maps'' $\Fi$ and $\hFi$ as
\be\label{PFi}
P^{\rm Ab}_n(p)=\li\Fi^n(\pi_p)\re
\quand
P^{\rm aB}_n(p)=1-\li\hFi^n(\pi_{1-p})\re
\qquad\big(p\in[0,1],\ n\in\N\big),
\ee
where the first equality follows from (\ref{Pcell}), Lemma~\ref{L:column}, and the fact that $(X_t)_{t\geq 0}$ started in product measure with intensity $p$ is the cellular automaton $(X^{{\rm Ab},p}_t)_{t\geq 0}$ from Subsection~\ref{S:cell}, and in the second equality we have used (\ref{hataB}).

Recall that $S=\{0,1\}^\N$ is equipped with the product topology, $\Mi$ is equipped with the topology of weak convergence, and $\Si\sub\Mi$ denotes the set of stationary measures.

\bl[Properties of the column maps]
The\label{L:Fiprop} maps $\Fi_{\rm A}$, $\Fi_{\rm b}$, $\Fi$, and $\hFi$ are continuous with respect to the topology on $\Mi$ and map $\Si$ into itself.
\el

\bpro
In view of (\ref{Fi}) it suffices to prove the statement for $\Fi_{\rm A}$ and $\Fi_{\rm b}$. Preservation of stationarity is obvious. To prove continuity, assume that $\mu_n\Rightarrow\mu$. By Skorohod's representation theorem \cite[Thm~6.7]{Bil99}, we can find random variables $Y_n,Y$ with laws $\mu_n,\mu$ such that $Y_n\to Y$ a.s.\ with respect to the topology on $S$. Let $(Y'_n)_{n\in\N}$ be an independent copy of $(Y_n)_{n\in\N}$ and let $Y'$ be the a.s.\ limit of $Y'_n$. Then $\Phi_{\rm A}(Y_n,Y'_n)\to\Phi_{\rm A}(Y,Y')$ a.s.\ and $\Phi_{\rm b}(Y_n)\to\Phi_{\rm b}(Y)$ a.s.\ with respect to the topology on $S$, which implies that $\Fi_{\rm A}(\mu_n)\Rightarrow\Fi_{\rm A}(\mu)$ and $\Fi_{\rm b}(\mu_n)\Rightarrow\Fi_{\rm b}(\mu)$.
\epro

\bl[Upper bound on the density]
One\label{L:nul} has
\be\label{FAbbd}
\li\Fi_{\rm A}(\mu)\re=\li\mu\re^2\quand\li\Fi_{\rm b}(\mu)\re\leq 2\li\mu\re\qquad(\mu\in\Si),
\ee
and consequently,
\be
\li\Fi(\mu)\re\leq 4\li\mu\re^2\quand\li\hFi(\mu)\re\leq 2\li\mu\re^2\qquad(\mu\in\Si).
\ee
\el

\bpro
Let $Y$ and $Z$ be independent random variables with law $\mu\in\Si$. Then
\be
\P\big[\Phi_{\rm A}(Y,Z)(j)=1\big]=\P\big[Y(j)\wedge Z(j)=1\big]=\li\mu\re^2\qquad(j\in\N),
\ee
and
\be
\P\big[\Phi_{\rm b}(Y)(j)=1\big]=\P\big[Y(j)\vee Y(j+1)=1\big]\leq\P\big[Y(j)=1\big]+\P\big[Y(j+1)=1\big]=2\li\mu\re
\ee
$(j\in\N)$, proving (\ref{FAbbd}). It follows that
\be\ba{l}
\dis\li\Fi(\mu)\re=\li\Fi_{\rm A}\circ\Fi_{\rm b}(\mu)\re=\li\Fi_{\rm b}(\mu)\re^2\leq\big(2\li\mu\re\big)^2,\\[5pt]
\dis\li\hFi(\mu)\re=\li\Fi_{\rm b}\circ\Fi_{\rm A}(\mu)\re\leq 2\li\Fi_{\rm A}(\mu)\re=2\li\mu\re^2.
\ec
\epro

Let $\un 0,\un 1\in S$ denote the configurations that are constantly zero and one, respectively. Since $\Phi_{\rm A}(\un 0,\un 0)=\un 0$ and $\Phi_{\rm b}(\un 0)=\un 0$, the delta-measure $\de_{\un 0}$ on $\un 0$ is a fixed point of the maps $\Fi_{\rm A}$ and $\Fi_{\rm b}$ and hence also of $\Fi$ and $\hFi$. The same applies to the delta-measure $\de_{\un 1}$ on $\un 1$. We define the domains of attraction of these fixed points as
\be\ba{ll}
\dis\Di_0:=\big\{\mu\in\Si:\Fi^n(\mu)\Asto{n}\de_{\un 0}\big\},
\quad&\dis
\Di_1:=\big\{\mu\in\Si:\Fi^n(\mu)\Asto{n}\de_{\un 1}\big\},\\[5pt]
\dis\hDi_0:=\big\{\mu\in\Si:\hFi^n(\mu)\Asto{n}\de_{\un 0}\big\},
\quad&\dis
\hDi_1:=\big\{\mu\in\Si:\hFi^n(\mu)\Asto{n}\de_{\un 1}\big\},
\ec
where $\Rightarrow$ denotes weak convergence of probability measures on $S$. The following proposition is the main result of this subsection.

\bp[Domains of attraction]
$\Di_0$\label{P:domain} and $\hDi_0$ are open, dense subsets of~$\Si$.
\ep

\bpro
Let $\Oi:=\{\mu\in\Si:\li\mu\re<1/4\}$. Since $S\ni x\mapsto x(0)$ is a bounded continuous function, the map $\Si\ni\mu\mapsto\li\mu\re:=\int_Sx(0)\mu(\di x)$ is continuous with respect to weak convergence of probability measures, and hence $\Oi$, being the inverse image of $[0,1/4)$ under this map, is an open subset of $\Si$. We claim that
\be
\Di_0=\big\{\mu\in\Si:\Fi^n(\mu)\in\Oi\mbox{ for some }n\geq 1\big\}.
\ee
Indeed, the inclusion $\sub$ follows from the fact that $\Oi$ is an open set containing $\de_{\un 0}$ and the inclusion $\supset$ follows from the fact that $\Oi\sub\Di_0$ by Lemma~\ref{L:nul}. Since the map $\Fi$ is continuous by Lemma~\ref{L:Fiprop}, the inverse image $\Fi^{-n}(\Oi)$ of $\Oi$ under the map $\Fi^n$ is open for each $n$, and hence
\be
\Di_0=\bigcup_{n\geq 1}\Fi^{-n}(\Oi),
\ee
being a union of open sets, is open. The same argument shows that also $\hDi_0$ is an open subset of $\Si$.

We observe that the map $\Fi_{\rm b}$ is linear but $\Fi_{\rm A}$ is not. It is easy to check that for any $\mu\in\Si$, we have
\be
\Fi_{\rm A}\big(r\mu+(1-r)\pi_0\big)=r^2\Fi_{\rm A}(\mu)+(1-r^2)\pi_0.
\ee
Using the linearity of $\Fi_{\rm b}$, we see that the same formula holds with $\Fi_{\rm A}$ replaced by $\Fi$ or $\hFi$. It follows that
\be
\Fi^n\big(r\mu+(1-r)\pi_0\big)=r^{n+1}\hFi(\mu)+(1-r^{n+1})\pi_0\Asto{n}\pi_0\quad\forall r\in[0,1),
\ee
which shows that $r\mu+(1-r)\pi_0\in\Di_0$ for all $r\in[0,1)$. Since $r\mu+(1-r)\pi_0\Rightarrow\mu$ as $r\to 1$, and $\mu\in\Si$ is arbitrary, we conclude that $\Di_0$ is dense in $\Mi$. The proof for $\hDi_0$ is the same.
\epro

\subsection{Proof of the main results}

In this subsection we prove the results stated in Subsection~\ref{S:intro} that have not been proved yet. Proposition~\ref{P:treecomp} has already been proved in Subsection~\ref{S:treecomp} and Proposition~\ref{P:ab} follows from Propositions \ref{P:APei} and \ref{P:BPei}, so it remains to prove Theorem~\ref{T:threshold} and Proposition~\ref{P:bounds}.\med

\bpro[of Theorem~\ref{T:threshold} and Proposition~\ref{P:bounds}]
We first prove the results for the game ${\rm Ab}_n(p)$. Since the map $L^{\rm Ab}_{2n}$ from (\ref{FAB}) is monotone with $L^{\rm Ab}_{2n}(\un 0)=0$ and $L^{\rm Ab}_{2n}(\un 1)=1$, formula (\ref{Pminmax}) tells us that the set
\be
I:=\big\{p\in[0,1]:P^{\rm Ab}_n(p)\asto{n}0\big\}
\ee
is an interval containing $0$ but not $1$. By formula (\ref{PFi}),
\be
I=\big\{p\in[0,1]:\pi_p\in\Di_0\big\},
\ee
so Proposition~\ref{P:domain} tells us that $I$ is an open set. It follows that there exists a constant $0<p^{\rm Ab}_{\rm c}\leq 1$ such that
\be\label{pab}
P^{\rm Ab}_n(p)\asto{n}0\quad\mbox{if and only if }p<p^{\rm Ab}_{\rm c}.
\ee
We recall from the discussion above Theorem~\ref{T:sharp} that $P^{\rm Ab}_n$ is strictly increasing and so also invertible. We observe that if $P^{\rm Ab}_n(p)=\ffrac{1}{4}$, then by formula (\ref{PFi}) and Lemma~\ref{L:nul},
\be
P^{\rm Ab}_{n+1}(p)=\li\Fi^{n+1}(\pi_p)\re\leq 4\li\Fi^n(\pi_p)\re^2=\ffrac{1}{4}\qquad(n\geq 1),
\ee
which proves that
\be
(P^{\rm Ab}_n)^{-1}(\ffrac{1}{4})\leq(P^{\rm Ab}_{n+1})^{-1}(\ffrac{1}{4})\qquad(n\geq 1).
\ee
It follows that the increasing limit
\be
p'_{\rm c}:=\lim_{n\to\infty}(P^{\rm Ab}_n)^{-1}(\ffrac{1}{4})
\ee
exists in $[0,1]$, and using Theorem~\ref{T:sharp} we see that
\be
(P^{\rm Ab}_n)^{-1}(q)\asto{n}p'_{\rm c}\qquad\forall 0<q<1,
\ee
which implies that
\be
P^{\rm Ab}_n(p)\asto{n}\left\{\ba{ll}
0\quad&\mbox{if }p<p'_{\rm c},\\
1\quad&\mbox{if }p>p'_{\rm c}.
\ea\right.
\ee
Combining this with (\ref{pab}) we see that $p'_{\rm c}=p^{\rm Ab}_{\rm c}$. Proposition~\ref{P:APei} implies that $p^{\rm Ab}_{\rm c}\leq 7/8$, while the bound $1/2\leq p^{\rm Ab}_{\rm c}$ follows from the simple argument given in the text below Proposition~\ref{P:bounds}. This completes the proofs for the game ${\rm Ab}_n(p)$.

For the game ${\rm aB}_n(p)$, the arguments are similar. In this case, by formula (\ref{PFi}),
\be
I:=\big\{p\in[0,1]:P^{\rm aB}_n(p)\asto{n}1\big\}=\big\{p\in[0,1]:\pi_{1-p}\in\hDi_0\big\},
\ee
is an open interval containing $1$ but not $0$. Combining this with Lemma~\ref{L:nul} and Theorem~\ref{T:sharp}, the arguments above show that there exists a $0\leq p^{\rm aB}_{\rm c}<1$ such that
\be
P^{\rm aB}_n(p)\asto{n}\left\{\ba{ll}
0\quad&\mbox{if }p<p^{\rm aB}_{\rm c},\\
1\quad&\mbox{if and only if }p>p^{\rm aB}_{\rm c}.
\ea\right.
\ee
Proposition~\ref{P:BPei} implies that $1/16\leq p^{\rm aB}_{\rm c}$ and the bound $p^{\rm aB}_{\rm c}\leq\ha(3-\sqrt{5})$ follows from (\ref{Pearl}) and Proposition~\ref{P:treecomp}.
\epro

\subsubsection*{Acknowledgments}

N.~Cardona-Tob\'on would like to express her gratitude to the University of G\"ottingen, where she was a postdoctoral researcher during the preparation of this paper. Additionally, she would like to thank the Institute of Information Theory and Automation for its hospitality during her research stay in 2023, which was supported by the Dorothea Schl\"ozer program of the University of G\"ottingen. J.M.~Swart is supported by grant 22-12790S of the Czech Science Foundation (GA CR). We thank Gerold Alsmeyer for bringing the reference \cite{ADN05} to our attention.

\end{document}